\documentclass[onefignum,onetabnum]{siamart220329}

\usepackage{xcolor}
\usepackage{amsfonts}
\usepackage{graphicx}
\usepackage{epstopdf}
\usepackage{overpic}
\usepackage{amsmath,esint}
\usepackage{bm}
\usepackage{epsfig}
\usepackage{epstopdf}
\usepackage{todonotes}
\usepackage{tikz}
\usetikzlibrary{arrows}
\usetikzlibrary{fadings}
\usepackage{pgfplots}
\pgfplotsset{compat=1.10}
\usepgfplotslibrary{fillbetween}
\usetikzlibrary{patterns}
\usepackage{enumerate}
\usepackage{enumitem}
\usepackage[mathscr]{euscript}
\usepackage{algorithm, algpseudocode}
\usepackage{mathtools}
\usepackage{appendix}
\renewcommand{\vec}[1]{\mathbf{#1}}
\newsiamremark{remark}{Remark}

\title{A time-frequency method for acoustic scattering with trapping} 
\author{ Heather Wilber\thanks{Department of Applied Mathematics, University of Washington (\email{hdw27@uw.edu}). Supported by the National Science Foundation under Grants DMS-2410045, DMS-2103317.} \and Wietse Vaes\thanks{ Department of Applied Mathematics, University of Washington (\email{wietsev@uw.edu}). Supported by the National Science Foundation under Grant DMS-2410045.} \and
Abinand Gopal\thanks{Department of Mathematics, University of California, Davis (\email{gopal@ucdavis.edu}).}
\and Gunnar Martinsson\thanks{Department of Mathematics and Oden Institute, University of Texas at Austin (\email{pgm@oden.utexas.edu}).} 
Supported by the Office of Naval Research (N00014-18-1-2354), by the National Science Foundation (DMS-1952735 and DMS-2313434), and by the Department of Energy ASCR (DE-SC0022251).} 

\begin{document}
\maketitle

\begin{abstract}
A Fourier transform method is introduced for a class of hybrid time-frequency methods that solve the acoustic scattering problem in regimes where the solution exhibits both highly oscillatory behavior and slow decay in time.  This extends the applicability of hybrid time-frequency schemes to domains with trapping regions. A fast sinc transform technique for managing highly oscillatory behavior and long time horizons is combined with a contour integration scheme that improves smoothness properties in the integrand. 
\end{abstract}



\section{Introduction} 

The acoustic scattering problem asks for the computation of a scattered wavefield that is produced when an incoming incident wavefield strikes a scattering object. Higher-order methods for these and related waveform simulations~\cite{maggio1994acoustic} have applications in areas such as virtual acoustic reconstruction~\cite{lumini2023project}. When sound--soft boundary conditions are imposed on the total field $u_{tot} = u_{inc} + u$, the scattered field $u$ satisfies the following initial boundary value problem: 
\begin{align} 
\label{eq:IVBP} 
& \dfrac{\partial^2 u}{\partial t^2}(x, t) - c^2\Delta u(x,t) = 0, \quad & (x, t) \in \Omega \times [0, T],\\
& u(x,0) = \dfrac{\partial u}{\partial t}(x,0) = 0, \quad &x \in \Omega, \label{eq:IVBP2} \\ \label{eq:IBVP3}
& u(x, t) = -u_{inc}(x,t), \quad & (x, t) \in \partial \Omega \times [0, T]. 
\end{align}
Here, $u_{inc}$ is the incident field, $c$ is the wave speed associated with the exterior domain $\Omega$, and $[0, T]$ is some relevant time interval over which the scattered wavefield is numerically detectable.  This problem is challenging to solve for a number of reasons. Standard domain discretization schemes (e.g., finite elements, finite volumes/differences) must manage the fact that the exterior domain is unbounded and impose artificial boundary conditions that mimic the behavior of the original problem (e.g.~the absorbing boundary layers method of~\cite{engquist1977absorbing}).
Traditional time-stepping methods must combat pervasive issues related to accumulating dispersive error and potentially prohibitive CFL restrictions on time step sizes. These issues are severely exacerbated when the domain includes corners that require intensive spatial refinement, and/or cavities where scattered waves decay in magnitude very slowly. 

Under the condition that the incident wavefield is enveloped by a Gaussian or is otherwise approximately bandlimited,  some of these complications  can be avoided by using \textit{hybrid time-frequency solvers}~\cite{anderson2020high, mecocci2000new}. In this class of solvers, the solution $u(x,t)$ is represented in terms of its Fourier transform $\hat{U}(x,\omega)$, which can be evaluated directly across a band $[W_1, W_2]$ of relevant frequencies by solving the boundary value problem
\begin{align} 
& \Delta \hat{U}(x, \omega) +\frac{\omega^2}{c^2} \hat{U}(x, \omega) = 0, \quad x \in \Omega, \label{eq:HH} \\
& \hat{U}(x, \omega) = -\hat{U}_{inc}(x, \omega), \quad x \in \partial \Omega,\label{eq:boundarycond} \
\end{align}
where $\hat{U}_{inc}(x, \omega)$ is the Fourier transform of $u_{inc}(x,t)$ with respect to $t$. The Sommerfeld radiation condition~\cite{schot1992eighty} is additionally imposed to exclude incoming waves from infinity.  The numerical task becomes the evaluation of the transform integral 
\begin{equation}
\label{eq:iFT}
 u(x,t) = \frac{1}{2 \pi} \int_{-\infty}^{\infty} \hat{U}(x, \omega) e^{-i \omega t} \mathrm{d} \omega \approx \frac{1}{2 \pi} \int_{W_1}^{W_2} \hat{U}(x, \omega) e^{-i \omega t} \mathrm{d} \omega, 
\end{equation}
where $\hat{U}(x,\omega)$ is the solution to the Helmholtz equation at frequency $\omega$ and point $x \in \Omega$. Here we assume that $W_1 > 0$ is bounded away from the origin to avoid a weak singularity in $\hat{U}(x, \omega)$ as $\omega \to 0$. However, there are specialized
quadrature methods that can be applied in such cases~\cite{Anderson_1984}. For completeness, we note that the other half of the transform pair is

\begin{equation}
\label{eq:FT}
\hat{U}(x, \omega) = \int_{-\infty}^{\infty} u(x,t) e^{i \omega t} \mathrm{d} t,
\end{equation}
where we define $u(x,t) = 0$ when $t<0$. 
\begin{figure}[t!]
 \centering
 \includegraphics[scale = .18, trim={3cm 1cm 2cm 0cm},clip]{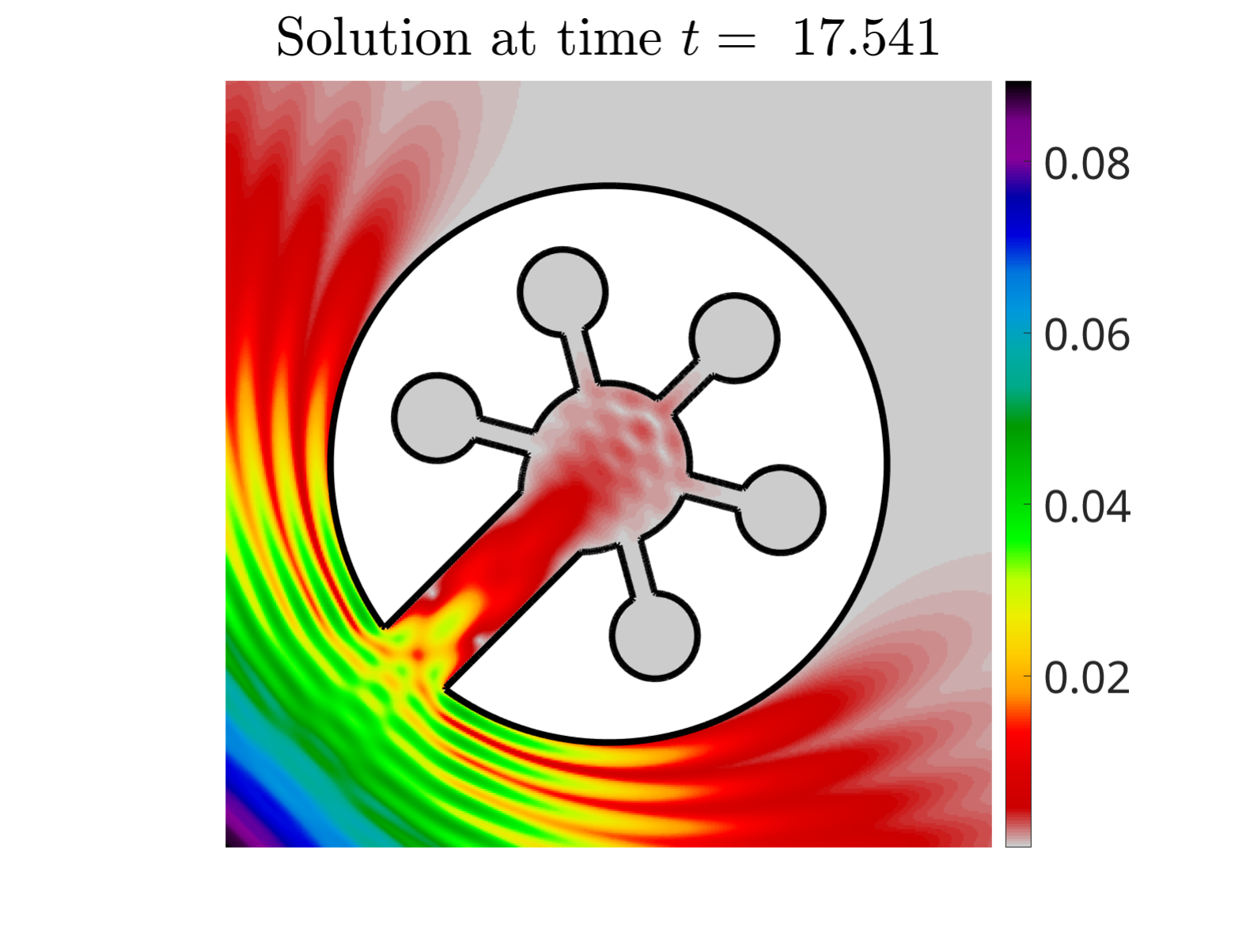}
  \includegraphics[scale = .18, trim={3cm 1cm 2cm 0cm},clip=true]{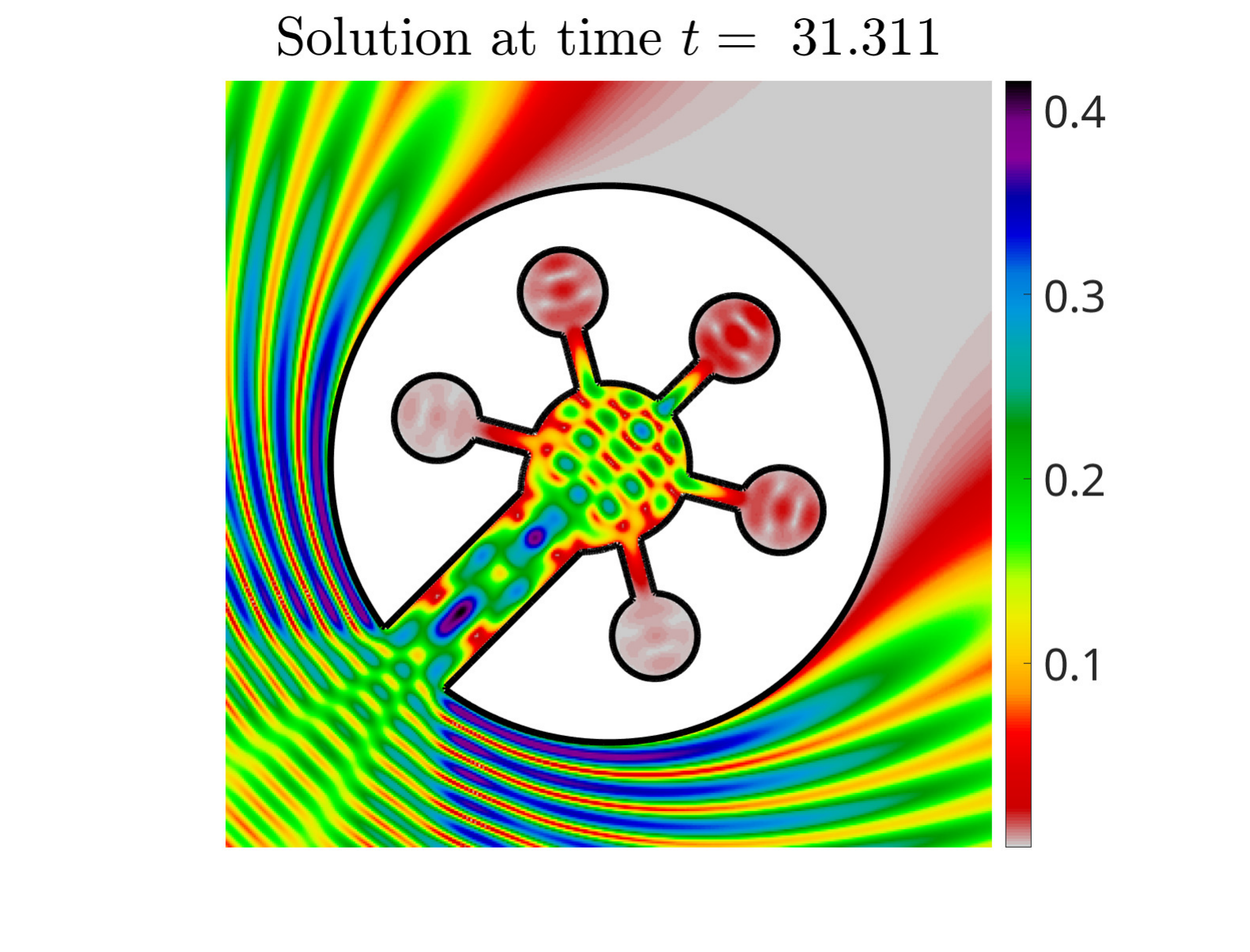}
   \includegraphics[scale = .18, trim={3cm 1cm 2cm 0cm},clip=true]{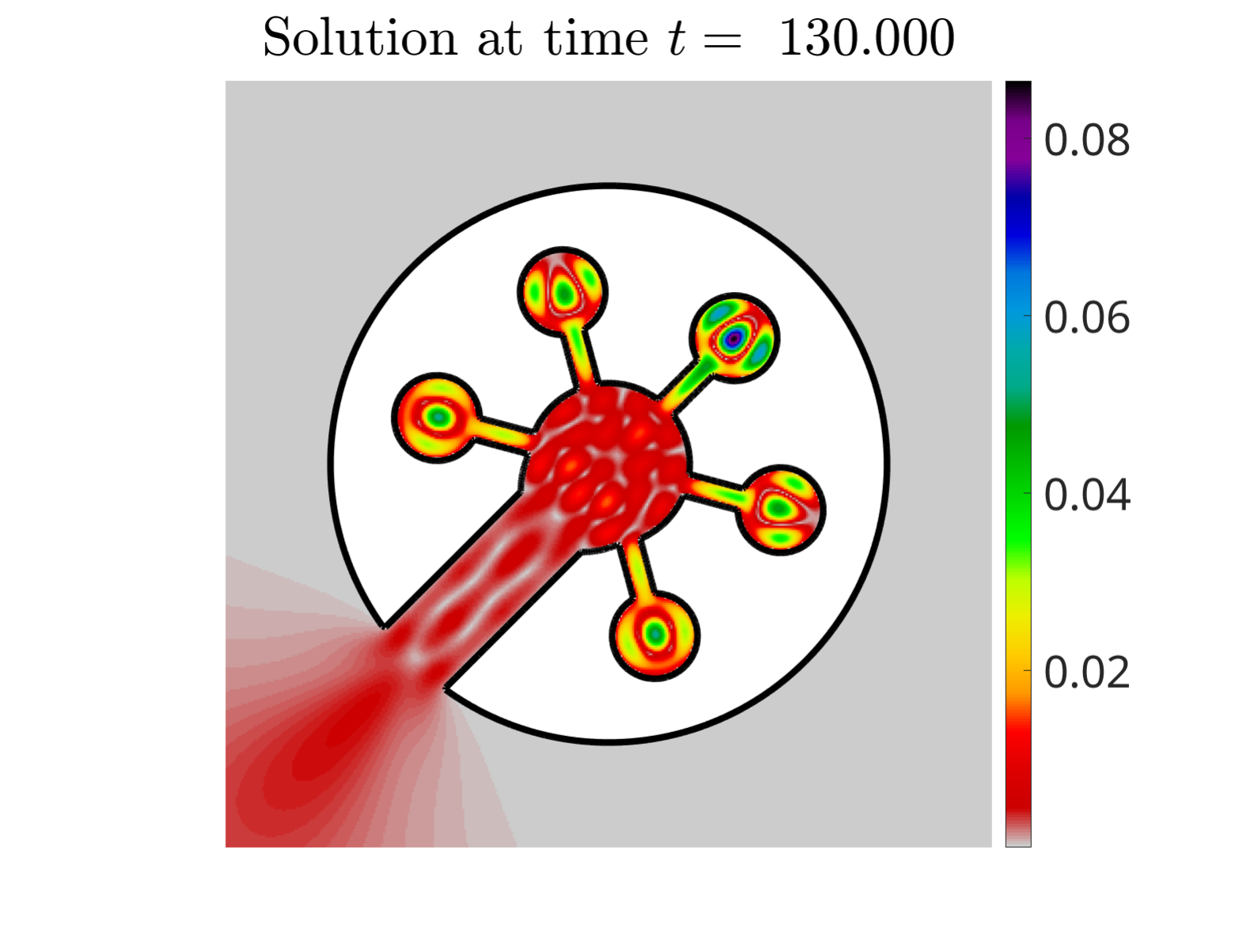}
\caption{ Plots of $|u_{tot}(x,t)|$ at times $t = 17.541$ (left), $t = 31.311$ (center), and $t = 130.00$ (right) on a ``magnetron" domain.  }
  \label{fig:gasket_magplots}
 \end{figure}

The time-frequency approach enjoys several advantages over more traditional time-domain-based methods, including the avoidance of restrictive CFL conditions, the ability to evaluate the solution at arbitrary points in time (time-skipping), and the production of solutions that are virtually free of dispersive error~\cite{anderson2020high}. Moreover, issues related to discretizing the exterior domain are eliminated since solutions to the Helmholtz equation can be expressed as integrals over $\partial \Omega$. Substantial progress in the development of practical time-frequency methods has been made in~\cite{anderson2020high}; related earlier work exists both in the context of Fourier transforms~\cite{mecocci2000new} and Laplace transforms~\cite{ banjai2012wave, lubich1994multistep}. We will focus on the Fourier-based methods, which can become prohibitively expensive unless strong assumptions about the smoothness of the domain and analyticity properties of $\hat{U}(x,\omega)$ are satisfied. They struggle in domains with corners and trapping regions, which are precisely the settings where the benefits of time-frequency hybrid methods are most needed. 

We introduce developments that make hybrid time-frequency methods effective in domains with these  ``unfriendly" features. To handle trapping, we develop a damping+correction scheme that combines contour deformation with fast transforms (the FFT and the fast sinc transform) to make evaluation of the integral in~\eqref{eq:iFT} tractable. This makes it possible to handle scenarios where solving the Helmholtz equation is more computationally intensive, such as in the ``magnetron" domain shown in \Cref{fig:gasket_magplots}.   Example domains where our method works well include domains with multiple scatterers (e.g., whisper galleries or acoustic panels), as well as domains with keyhole regions, channels, and multiple corners or cusps. 

We focus on analyzing and resolving two central problems related to trapping: 

\vspace{.1cm}
 \begin{itemize}
   \item ({\bf P1}) As $t$ grows, the integral in~\eqref{eq:iFT} becomes highly oscillatory. Standard quadrature rules (e.g., trapezoidal or Gauss-Legendre) require a number of quadrature points that must grow linearly with $t$ to achieve any accuracy. 
   \item ({\bf P2})   $\hat{U}(x, \omega)$ may have poor analyticity properties, with nearby singularities in the lower half of the complex plane. This can make the convergence of quadrature rules applied to~\eqref{eq:iFT} very slow. 
 \end{itemize} 
 \vspace{.1cm}

\begin{figure}[h!]
    \begin{minipage}{.72\textwidth} 
 \centering
  \begin{overpic}[width=\textwidth]{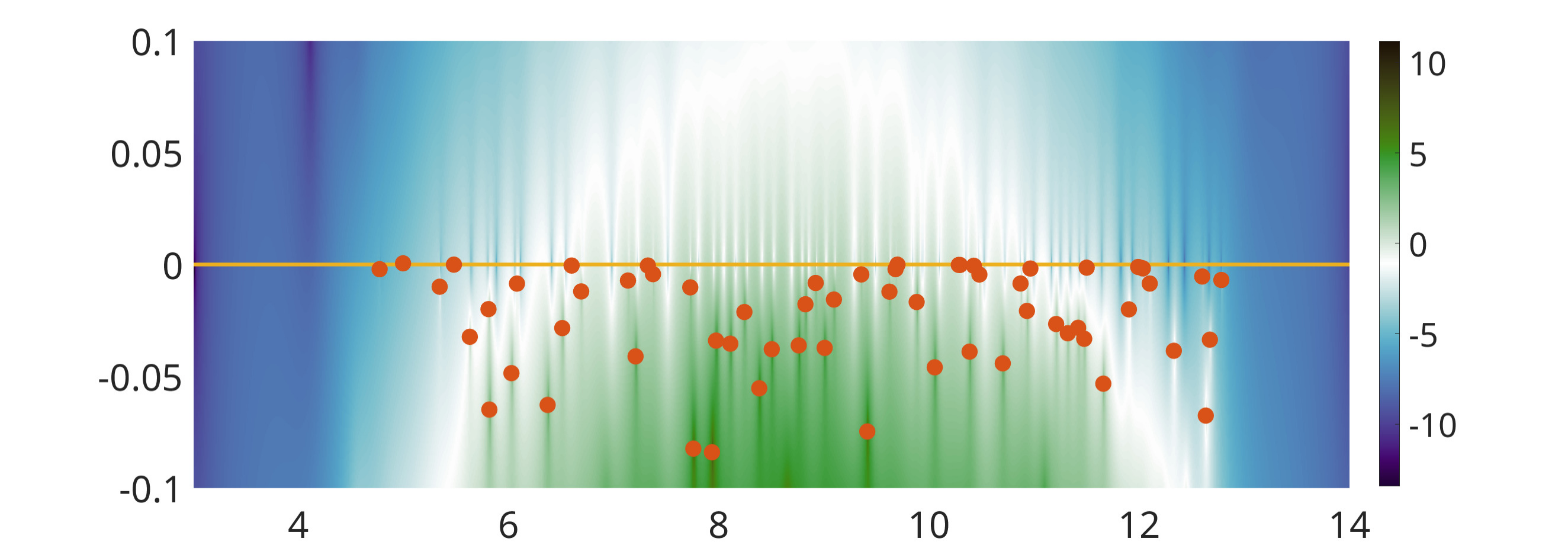}
  \put(45, -3){ {\small $Re(\omega)$ }}
  \put(1, 10){\rotatebox{90}{{ \small $Im(\omega)$}}}
\end{overpic}
	\end{minipage}
   \begin{minipage}{.27\textwidth} 
\includegraphics[scale = .1, trim={0cm 1cm 4cm 1cm},clip]{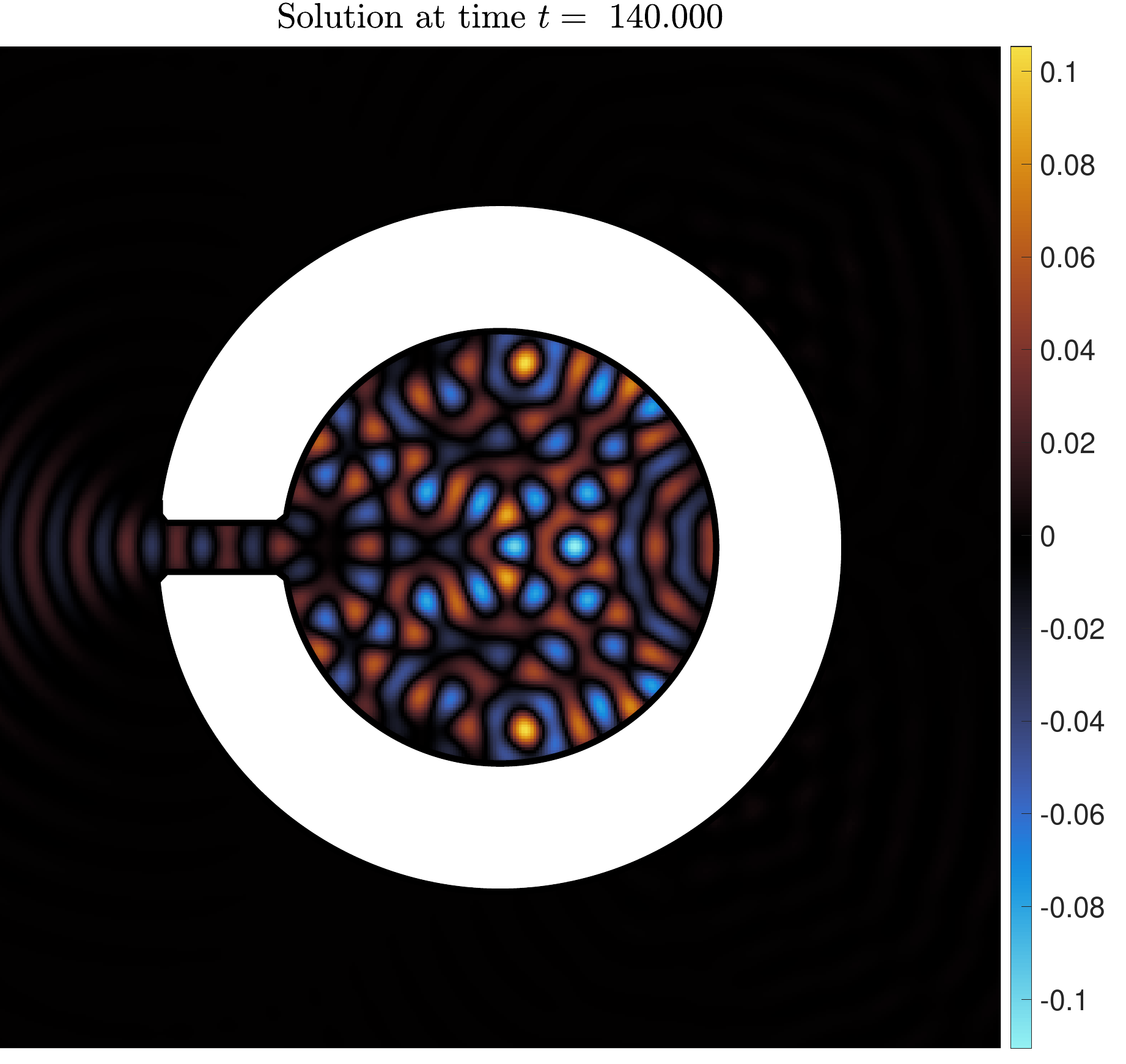}
	\end{minipage}
\caption{ \textit{Left:} A plot of the magnitude of a meromorphic extension of $\hat{U}(\tilde{x}, \omega)$, associated with the domain on the right, is shown in the complex $\omega$-plane around the frequency interval $[3, 14]$ for a single point $\tilde{x} = [0,0]^T$. The colors are displayed using a logarithmic scale (e.g., from $10^{-11}$ (deep blue) to $10^{10}$ (dark green)). The point $\tilde{x} = [0, 0]^T$ is center of the keyhole cavity depicted on the right. Poles of the extended function, which correspond to singularities of $\hat{U}(\tilde{x}, \omega)$, are depicted as orange dots. The imaginary part of the pole closest to the real line is $  \approx 1 \times 10^{-5}$ in magnitude.  \textit{Right:} The domain where $\hat{U}$ was computed is shown for reference, along with the real part of $u_{tot}(x,t)$ at $t = 68.121$. Details on the experiment can be found in \Cref{sec:num}.}
  \label{fig:poleplot}
 \end{figure}

Note that ({\bf P1}) and ({\bf P2}) are coupled: If $\tilde{x}$ is a fixed point in an activated trapping region of $\Omega$, then the decay of $|u(\tilde{x},t)|$ is slow, so we expect that large values of $t$ will be of interest. However, this is also precisely when one must worry about (\textbf{P2}).  As shown in \Cref{fig:poleplot}, poles of the meromorphic extension\footnote{computed via the AAA algorithm~\cite{nakatsukasa2018aaa}.} of $\hat{U}(\tilde{x}, \omega)$ in the complex $\omega$-plane creep closer to the real line as the trapping becomes more severe.  The distance of the nearest pole to the real line corresponds to the decay rate of $|u(\tilde{x}, t)|$; this phenomenon is connected to the analysis and computation of the eigenvalues and near-resonant modes of scattering operators~\cite{bruno2024evaluation, patel2023wigner}. 

To improve the situation, we note that singularities of $\hat{U}(x, \omega)$ always lie below the real line. This suggests the consideration of the perturbed function $\hat{U}(x, \omega+i\delta)$, where $\delta > 0$.  One can view this perturbation as the introduction of damping into the solution. Integrating over a rectangular contour (traversed counterclockwise) in the upper half plane, Cauchy's integral theorem implies that the solution can be represented as 
\begin{equation} \label{eq:rect}
\int_{W_1}^{W_2} \hat{U}(x, \omega) e^{-i \omega t} \mathrm{d} \omega  =  \underbrace{\int_{W_1}^{W_2}  \hat{U}(x, \omega + \delta i) e^{-i (\omega + \delta i) t} {\rm d} \omega}_{I_{\delta}} - I_{cL}-  I_{cR},
\end{equation}
where the correction terms $I_{cL}$ and $I_{cR}$ are integrals along the vertical sides of the rectangle.  There is an inherent upper limit on $\delta$. If it is taken too large,  the correction terms cannot be stably evaluated. However, the room afforded by $\delta$ is generous enough that when paired with a fast and accurate evaluation scheme for the highly oscillatory integral $I_\delta$, we are able to solve the acoustic scattering problem, even with severe trapping, over long time horizons. Combined with highly effective and well-established methods~\cite{gillman2012direct, martinsson2019fast} for solving the Helmholtz equation on challenging domains (e.g., with corners, cusps, etc.), this approach can manage complicated domains. This is illustrated with several examples in \Cref{sec:num}; videos of various simulations can be seen at the website listed in~\cite{HOMEyoutube}. 

 Our method can be applied using any (damped) Helmholtz solver that can handle complex-valued frequencies. To manage corners, one can incorporate state-of-the-art methods based on specialized series expansions~\cite{bruno2009high,hoskins2019numerical, serkh2019solution}, QBX-based methods~\cite{barnett2014evaluation} or more naive but easier-to-implement refinement schemes~\cite[Ch.~12]{martinsson2019fast}, which is what we use in our examples. One could also consider more sophisticated strategies, such as broadband Helmholtz solvers~\cite{gopal2022broadband}; we investigate this approach in greater depth in future work. While we emphasize 2D problems in this paper, the ideas extend naturally to the 3D setting (see \Cref{fig:cruellerplots}). For readability, it should be assumed unless otherwise stated that $\Omega$ is a subset of $\mathbb{R}^2$. 

 \subsection{Related work}

 Two papers were foundational in the development of our work. Rokhlin's 1983 paper~\cite{rokhlin1983solution} on solving acoustic scattering problems introduces ideas fundamental to the numerical development of boundary-integral-based approaches to this problem, and also applies perturbation to the frequency parameter to alleviate trapping. Rokhlin credits the latter idea to Courant and Hilbert~\cite{courant1962methods}.  Several other important numerical developments for representing and computing solutions (without strong trapping) with the Fourier-based time-frequency method were introduced in Anderson, Bruno and Lyon's 2020 paper~\cite{anderson2020high}, including a method for efficiently handling  highly oscillatory inverse Fourier integrals via truncated sinc expansions and fast transforms. The method we apply is mathematically equivalent to this one but considers the involved ideas from the perspective of classical results on sampling schemes for the recovery of bandlimited functions.
 
A recently published preprint from Bruno and Santana~\cite{bruno2025efficient} is the work closest to our own, as it tackles the same problem.  Here, the authors use a polynomial + rational representation of $\hat{U}(x, \cdot)$ to formulate a singularity subtraction method~\cite{bruno2025efficient} for handling issues with trapping. In contrast, our method uses a purely polynomial representation and applies damping to improve the generally slow convergence rate of polynomial approximation in this setting. The method in~\cite{bruno2025efficient} is shown to be effective for a variety of challenging problems. Moreover, it has the advantage that it directly identifies near-resonances and therefore allows one to derive large-time asymptotic expansions of the scattered field. However, the approach is based on AAA approximation, and the analysis of robust convergence with respect to the number of Helmholtz solves is delicate.  The ideas in~\cite{bruno2025efficient} can likely be combined with the ideas developed here. For instance, if the resolvent-based methods in~\cite{bruno2024evaluation} for pole-finding are adapted to handle complex-valued frequencies, then the damping + correction idea can be used in concert with a polynomial + rational representation of $\hat{U}(x, \cdot)$. This would reduce both the number of poles required to accurately represent the slow-decaying components of $u(x,t)$ and the number of evaluations of $\hat{U}(x, \cdot)$ required for constructing the rational component of the approximation.

Other less-directly related works include the hybrid Fourier approach from~\cite{mecocci2000new} and convolutional quadrature methods based on Laplace transforms~\cite{banjai2010multistep,  banjai2012wave, banjai2014fast, lubich1994multistep}. Approaches related to retarded potential integral functions for a boundary element method in the time domain can be found in~\cite{barnett2020high, et1986formulation, yilmaz2004time}.

 \subsection{Organization} The rest of this manuscript is organized as follows: In \Cref{sec:HTF}, we briefly review existing work on Fourier-based hybrid time-frequency methods, as well as relevant details related to solving the Helmholtz equation with boundary integral methods. \Cref{sec:fast-sinc} discusses the fast sinc-based quadrature scheme. \Cref{sec:DC} describes a damping + correction method for regions with trapping, and \Cref{sec:num} gives numerical results and examples.

\section{The Fourier-based hybrid time-frequency method} 
\label{sec:HTF}

After briefly discussing the evaluation of $\hat{U}_{inc}$, which is not the main focus of our work, we describe the Fourier-based hybrid time-frequency approach in an ideal setting where $|u(x,t)|$ decays rapidly and trapping is a nonissue. This gives us the opportunity to review fundamental aspects of the solver. Then in Sections \ref{sec:fast-sinc} and \ref{sec:DC}, we turn to domains with trapping.

\subsection{Transforming the incident wavefield} 
As seen in~\eqref{eq:boundarycond}, the incoming field $-\hat{U}_{inc}(x, \omega)$ furnishes Dirichlet boundary data.  The class of solvers we are interested in are applicable in regimes where $u_{inc}$ decays sufficiently rapidly in both frequency space and in the time domain. This means that at any fixed point $\tilde{x}$ in space, $u_{inc}(\tilde{x}, t)$ practically vanishes (e.g., \hbox{$|u_{inc}(\tilde{x},t)|< \epsilon_{mach}$}) outside a finite interval of support over\footnote{While the scattered field is zero-valued for $t < 0$, this is not the case for the incident field.}  $t \in (-\infty, \infty)$.  
Additionally, the Fourier transform 
  $\hat{U}_{inc}(\tilde{x}, \omega) = \mathcal{F}[u_{inc}(\tilde{x}, \cdot)](\omega)$, nearly vanishes outside a fixed band of frequencies $[W_1, W_2]$.  We take as a  model incident field the Gaussian packet of the form
  \begin{equation}
  \label{eq:enveloped_osc_inc_wave}
    g(x,t) = \frac{1}{\sqrt{2 \pi}\sigma} e^{ -\left( \frac{(t-s(x))^2}{2 \sigma^2}+ i \omega_0 (t - s(x) )\right)}, \quad s(x) = \frac{x \cdot z_0}{c} + t_0,
  \end{equation}
  where $w_0$ and $t_0$ are fixed  phase and time shift parameters chosen so that the initial conditions on $u$ in~\eqref{eq:IVBP2} are approximately satisfied over $x \in X$, and $z_0 \in \mathbb{R}^2$ is a fixed unit direction of propagation.  Throughout the paper, we use incident fields of the form~\eqref{eq:enveloped_osc_inc_wave}. In this case the function $\hat{U}_{inc}$ is known explicitly so evaluation is trivial.
  
  A broader class of incident fields can also be considered, but the transform integrals for computing the Fourier transform of these fields are potentially highly oscillatory. A major focus of~\cite{anderson2020high} is the development of a windowing method that alleviates this issue. The windowing method also results in a separation of scales in the frequency domain that can exploited in later stages of the hybrid time-frequency solver. These ideas can optionally be applied in our method, allowing for more versatility in the choice of incident waves. However, since this is not the focus of the present manuscript, we limit our consideration of incident wavefields. Given that $\hat{U}_{inc}$ is accessible, we consider the problem of evaluating $u(x,t)$ at all points $X \times S$, where $$X = \{x_1, x_2, \cdots, x_M\} \subset \Omega, \quad S = \{t_1, t_2, \cdots t_N\} \subset [0, T].$$

\subsection{The fast-decaying case} We first consider the setting where trapping is very mild or nonexistent. In such settings, the scattered waves dissipate quickly and there is no need to take $t$ large.

To numerically simulate the scattered wavefield in this scenario, a quadrature rule, such as Gauss-Legendre, is selected and the truncated integral in~\eqref{eq:iFT} is discretized so that for each $(j,k)$ pair,

\begin{equation} \label{eq:discSOL} u(x_j,t_k) \approx \frac{1}{2 \pi}  \sum_{\ell = 1}^m \eta_\ell \hat{U}(x_j, \omega_\ell) e^{-i \omega_\ell t_k}.\end{equation} Here, $\eta_\ell$ and $\omega_\ell$ are quadrature weights and nodes, respectively. Since $t_k$ is never especially large and $\hat{U}(x, \cdot)$ is sufficiently nice, the size of $m$ is constrained primarily by the frequency content in the problem. 
The central computational challenge in computing the sum in~\eqref{eq:discSOL} is the evaluation of $\hat{U}(x_j, \omega_\ell)$ over the set $\{\omega_\ell\}_{\ell = 1}^m \times \{ x_j\}_{j = 1}^M$, which requires solving the Helmholtz equation at each fixed frequency $\omega_\ell$. In principle, any Helmholtz solver can be used for this purpose. A wealth of such solvers exist that can deliver high accuracy solutions even in settings where weak singularities induced by corners arise. Examples include specialized solvers that use expansions of fundamental solutions~\cite{barnett2010exponentially, gopal2019}, as well as variations of fast direct solvers that use boundary integral formulations in conjunction with specialized discretization strategies for corner geometries~\cite{bruno2009high, serkh2019solution} and near-singular integrals~\cite{af2018adaptive, barnett2014evaluation,krantz2024error,  wu2021zeta}. We have implemented a solver using a boundary integral formulation, and give a brief overview of relevant details here. A more thorough review can be found in~\cite{martinsson2019fast}. 
\subsubsection{A boundary integral formulation}\label{sec:BIE}

Consider a single fixed frequency $\omega_\ell >0$ and set $\kappa_\ell = \omega_\ell/c$. The solution to~\eqref{eq:HH} at $\kappa_\ell$ over all $x \in \Omega$ can be expressed in terms of an integral equation that depends on an unknown density function \hbox{$\phi_\ell$~\cite[Ch.~10]{martinsson2019fast}}. With the so-called combined field formulation, the solution is 
\begin{equation} 
\label{eq:U_int}
\hat{U} (x, \omega_\ell)  = \int_{\partial \Omega} \left(d_{\kappa_\ell}(x,y) - i Re(\kappa_\ell) s_{\kappa_\ell}(x,y) \right) \phi_\ell(y) {\rm ds}(y), 
\end{equation}
where 
\begin{equation}
s_{\kappa_\ell}(x,y) = G_{\kappa_\ell}(x-y), \quad d_{\kappa_{\ell}}(x,y) = n(y) \cdot \nabla_x G_{\kappa_\ell}(x-y), \quad G_{\kappa_\ell}(x) = \frac{i}{4} H_0^{(1)}(\kappa_\ell|x|).
\end{equation}
Here, $H_0^{(1)}$ is the zeroth-order Hankel function of the first kind~\cite[Sec.~10.2]{olver2010nist}, $n(y)$ is the outward normal unit vector to $\partial \Omega$, and $G_{\kappa_\ell}$ is the free-space fundamental solution to the 2D Helmholtz equation. The work required to solve~\eqref{eq:HH} over $\{ x_j\}_{j = 1}^M \times \{\omega_\ell\}_{\ell = 1}^m$ can be accomplished in two stages. First, we construct a set of density functions $\Phi = \{ \phi_1, \phi_2, \cdots \phi_m\}$. Then, we use the solution formulation in~\eqref{eq:U_int} to evaluate each $\hat{U}(x, \omega_\ell)$ at every $x \in X$.

The unknown function $\phi_\ell$ satisfies the following second-kind Fredholm equation on $L^2(\partial \Omega)$\cite[Ch.~11]{martinsson2019fast}: 
\begin{equation}
\label{eq:mixedBIE}
\frac{1}{2} \phi_\ell(x) + \int_{\partial \Omega} \left(d_{\kappa_\ell}(x,y) - i \kappa_{\ell} s_{\kappa_\ell}(x,y)\right)\phi_\ell(y) {\rm ds}(y) = -\hat{U}_{inc}(x, \omega_\ell).
\end{equation}
We choose to work with the ``combined field formulation'' rather than those based on only the single or double layer potential because this ensures that the problem is well-posed (there are no spurious resonances).

A Nystr\"om discretization of~\eqref{eq:mixedBIE} of size $n$ leads to an $n \times n$ linear system that can be solved for the vector $\tilde{\phi}_\ell = \begin{bmatrix} \phi_\ell(y_1) & \phi_\ell(y_2) & \cdots & \phi_\ell(y_n) \end{bmatrix}^T$, which we take as a representation of $\phi_\ell$. Here, $n$ grows linearly with the wave number $\kappa_\ell$.
For smooth boundaries, we apply the Zeta correction quadrature scheme from~\cite{wu2021zeta}, which is designed to manage the singularity in the integrand of the integral operator in~\eqref{eq:U_int}.  Modifications for nonsmooth boundaries are discussed in \Cref{sec:cornerdoms}. In this way, we compute and store representations for the set of density functions. 

In the second stage, vectors $\{b_\ell(\vec{x}) = \hat{U}(\vec{x},\omega_\ell)\}_{\ell = 1}^m$ are computed, where \newline \hbox{$\vec{x} = \begin{bmatrix} x_1,& \cdots &x_M\end{bmatrix}^T$}, by evaluating~\eqref{eq:U_int}.  If the discretization of each contour integral  in~\eqref{eq:U_int} requires $\leq n$ quadrature points, then the naive evaluation of each integral requires a matrix vector product with an $M \times n$ matrix. The total computational cost for all $b_\ell(\vec{x})$ is $\mathcal{O}(Mnm)$. However, this cost can be reduced to only $\mathcal{O}(m(M +n))$ by instead using the fast multipole method (FMM)~\cite{greengard1987fast}. 

With all of the $ b_\ell(\vec{x}) \approx \hat{U}(\vec{x}, \omega_\ell)$ values computed, the solution $u(\vec{x}, t)$ in the time domain can then be evaluated at each $t_k$ via~\eqref{eq:discSOL}.  This requires the evaluation of $MN$ sums of the form
\begin{equation} \label{eq:dft3}
 \sum_{\ell = 1}^{m} \eta_\ell b_\ell(x_j) e^{- i \omega_\ell t_k}.
 \end{equation}
Naively this takes an additional $\mathcal{O}(NMm)$ computations. However, since~\eqref{eq:dft3} has the form of a type-III nonuniform discrete Fourier transform, we can make use of the NUFFT-III~\cite{barnett2018non} to do the computation in only $\mathcal{O}(M(P \log P +m + N))$ computations, where $P = W_2-W_1$. Now we turn to the setting where $|u(x,t)|$ decays slowly in time. 

\section{Fast sinc quadrature for large time evaluation} \label{sec:fast-sinc}

We first consider (\textbf{P1}) without regard for (\textbf{P2}). A basic approach is to consider the worst-case scenario where $t = T$, and then discretize the integral in~\eqref{eq:discSOL} with as many points as is needed. Due to the $e^{-i \omega t}$ term, the number of points required to achieve any accuracy must scale linearly with $T$. An alternative approach, which turns out to be equivalent to  the ``scaled convolution" method presented in~\cite{anderson2020high}, is readily suggested by Shannon–Whittaker–Kotelnikov Sampling Theorem~\cite{shannon1949communication, whittaker1915xviii}. We present the idea from the perspective of Shannon sampling, as this usefully links the decay behavior of $u(x, t)$ to smoothness and analyticity properties of $\hat{U}$.  

The sampling theorem states that in the Paley-Wiener space consisting of functions $f$ such that $f$ is continuous on $\mathbb{R}$ with bounded $L_2$ norm and \newline \hbox{${\rm supp}(\mathcal{F}f) \subseteq [-B, B]$}, $f$ can be expressed as

\begin{equation} \label{eq:sincresult}
    f(t) = \sum_{j \in \mathbb{Z}} f\left( \tfrac{j}{L} \right) {\rm sinc}(L t -j ) , \quad L \geq 2B,
\end{equation} 
with ${\rm sinc}(t) = \sin(\pi t) / (\pi t)$.  In general, truncations of~\eqref{eq:sincresult}  converge slowly due to the sluggish decay rate of the sinc function. However, in our settings the function samples $|u(x,t_k)|$ decay rapidly enough (see \Cref{fig:smoothing_slices}) to make approximations based on~\eqref{eq:sincresult} more reasonable. 

We cannot sample $u(x,t)$ directly since it is what we are trying to compute. To convert~\eqref{eq:sincresult} to a practical method, we note that if $\hat{U}(x, \omega)$ were truly bandlimited, then the Fourier coefficients of the periodic extension of $\hat{U}(x,\cdot)$ on $[W_1, W_2]$ can exactly be expressed using equally-spaced samples of $u(x, t)$ over $\mathbb{R}$. Letting $y \in [0, 1]$, we apply the change of variables $\omega = Py + W_1$, where $P = W_2-W_1$. There are coefficient functions $\{c_{-m}(x), \cdots, c_{m}(x)\}$ such that 
 \begin{equation} 
 \label{eq:trigpoly}
 \hat{U}(x, Py + W_1) \approx \frac{1}{2m + 1}\sum_{j = -m}^{m} c_j(x) e^{2 \pi  i j y },
 \end{equation}
 where $$ c_j(x) = \int_{0}^{1} \hat{U}(x, Py+W_1) e^{-2 \pi i j y} {\rm d}y. $$
The coefficients can be computed for any fixed $x$ by sampling $\hat{U}(x, \cdot)$ on $2m \!+\!1$ equally-spaced points over $[W_1, W_2)$ and then applying the FFT. In the truly compactly supported case, it follows from an inspection of~\eqref{eq:iFT} that
\begin{equation} \label{eq:candtime}
    c_j(x) = \frac{2 \pi}{P} e^{2 \pi i j W_1/P} (-1)^j u(x, 2 \pi j/P).
\end{equation}
From this we see that for $j \leq 0$, $c_j(x) = 0$. Substituting~\eqref{eq:trigpoly} into the integral in~\eqref{eq:iFT} yields the following approximation for $u$: 
 \begin{equation}   
 \label{eq:sincsum}
u_m(x,t) = \frac{P}{2 \pi (2m+1)} e^{-i t(P/2 + W_1)} \sum_{j = 1}^m (-1)^j c_j(x)   {\rm sinc}\left( \tfrac{Pt}{2 \pi} - j \right).
 \end{equation}
Once the coefficients are known, the fast sinc transform~\cite{greengard2007fast,lawrenceFS} can be used to efficiently evaluate $u_m(x,t)$ at $N$ points $\{t_1, \ldots, t_N\} \subset [0, T]$ at each fixed point $x_k \in X$. The computational cost is $\mathcal{O}(P \log P +m + N)$ operations, and the FINUFFT-based implementation~\cite{lawrenceFS} applies batching techniques that substantially improve the practical run time when~\eqref{eq:sincsum} must be evaluated for many values of $x$.

 The approximation in~\eqref{eq:sincsum} is equivalently derived by convolving $u(x,t)$ with a scaled sinc function. This is perhaps the  simplest windowing approximation strategy among many possibilities, including the use of Kaiser-Bessel functions, Gaussians, bump functions, or prolate-spheriodal wave functions. Alternative strategies could potentially improve the rate at which \eqref{eq:sincsum} converges to $u(x,t)$ as $m$ grows~\cite{kircheis2024fast}, but we leave this for future investigation.  We emphasize that the primary advantage of~\eqref{eq:sincsum} over basic quadrature choices is that the $t$-dependent oscillatory content in the integrand is reduced to a multiplicative factor (outside the integrand) that does not directly dictate the rate at which $\hat{U}(x, \cdot)$ must be sampled. This solves  \textbf{(P1)}. 

\subsection{Convergence of the truncated sinc expansion} 

We now turn to \textbf{(P2)}. For the sake of analysis, we assume in this section that $\hat{U}(x, \cdot)$ is truly compactly supported on $[W_1, W_2]$. The fundamental challenge in using~\eqref{eq:sincsum} is that $\hat{U}(x, \cdot)$ becomes difficult to approximate with trigonometric polynomials if the decay of $|u(x,t)|$ is very slow. This in turn is closely connected to the smoothness and analyticity properties of $\hat{U}(x, \cdot)$.  We succinctly state the obvious coupling between the decay rate of $|u(x,t)|$ and the convergence of $u_m \to u$ in the following lemma, which is a direct consequence of the Shannon–Whittaker–Kotelnikov Sampling Theorem and~\eqref{eq:candtime}. 

\begin{lemma}  For a fixed $\tilde{x}\in \Omega$, let $\hat{U}(\tilde{x}, \omega)$ be a solution to~\eqref{eq:HH} that is compactly supported on $[W_1, W_2]$. Let $u(\tilde{x}, t) = \mathcal{F}^{-1}(\hat{U}(\tilde{x}, \omega))$ and define $u_m(\tilde{x}, t)$ as in~\eqref{eq:sincsum}. Let $0 < \epsilon < 1$.  Suppose there is a monotonically decreasing and invertible function  $g(t)$ such that $ |u(x,t)| = \Theta(g(t))$.
Then, 
  $$m = o\left( g^{-1}(\epsilon ) \right) \implies \| u(\tilde{x}, \cdot)-u_m(\tilde{x}, \cdot) \|_{\infty, [0, \infty)} = \mathcal{O}(\epsilon)$$.
\label{lemma:tbnd}
 \end{lemma}

 \Cref{lemma:tbnd} implies that we can only guarantee that $\|u(\tilde{x}, \cdot)-u_m(\tilde{x}, \cdot)\|_{\infty, [0, \infty)} \leq \epsilon$ if $m$ is roughly on the order of a time point $t_\epsilon$, where $|u(x, t_\epsilon)| < \epsilon$ and decays steadily thereafter for all $t > t_\epsilon$. 
 For example, if $\mathbb{R}^d \setminus \Omega$ does not possess trapping regions (e.g., the scatterer forms a star domain\footnote{A set $\mathcal{D} \subset \mathbb{R}^d$ is called a \textit{ star domain} if there exists $p_0 \in \mathcal{D}$ such that for all $p \in \mathcal{D}$, the line segment joining $p_0$ and $p$ also lies in $\mathcal{D}$.}) and $d$ is odd\footnote{For even $d$ the situation is slightly more complicated~\cite{morawetz1975decay}.}, then the celebrated result of Morawetz, Ralston and Strauss~\cite{morawetz1977decay} proves that $g(t)$ is an exponential function. This implies that  $u_m \to u$ uniformly and pointwise at an exponential rate. In such a setting, we expect the sinc-based method to work extremely well. However, we also expect that $T$ will never be required to be especially large in such a case.

 \begin{remark} The claim made in~\cite{anderson2020high} is that by using~\eqref{eq:sincsum}, the cost of evaluating the scattered field no longer appears to depend directly on the magnitude of the largest numerically relevant time $T$.  However, Lemma~\ref{lemma:tbnd} suggests that the choice of $m$ required to achieve a specified accuracy $\epsilon$ implicitly depends linearly on some  $t_\epsilon \in [0, T]$. If we take $T$ to be $t_\epsilon$, then this appears to be a contradiction. However, a subtlety not to be overlooked is that if we choose $\epsilon$ such that $t_\epsilon < T$, we have $\mathcal{O}(\epsilon)$ error, even for $t \in [t_\epsilon, T]$. This is a stark contrast with, say, Gauss-Legendre quadrature, where the error only decays meaningfully if the number of quadrature points used to discretize~\eqref{eq:iFT} is $\mathcal{O}(t)$, and thus in the worst case is $\mathcal{O}(T)$. It is in this sense that, as the authors of~\cite{anderson2020high} state, the accuracy of $u_m(x, t)$ over $[0, T]$ is \textit{tunable}.  

\end{remark}

What can we say about $g(t)$ in cases where there is trapping? The decay of $|u(x,t)|$ may be arbitrarily slow, and it cannot be succinctly described, though some asymptotics related to the poles of the scattering operator can be numerically extracted~\cite{bruno2025efficient, bruno2024evaluation}.
We observe that the convergence of $u_m \to u$ in this setting can be so slow that~\eqref{eq:sincsum} is no longer numerically useful. 

\begin{figure}
\centering
\hspace{1cm}
\includegraphics[scale = .165]{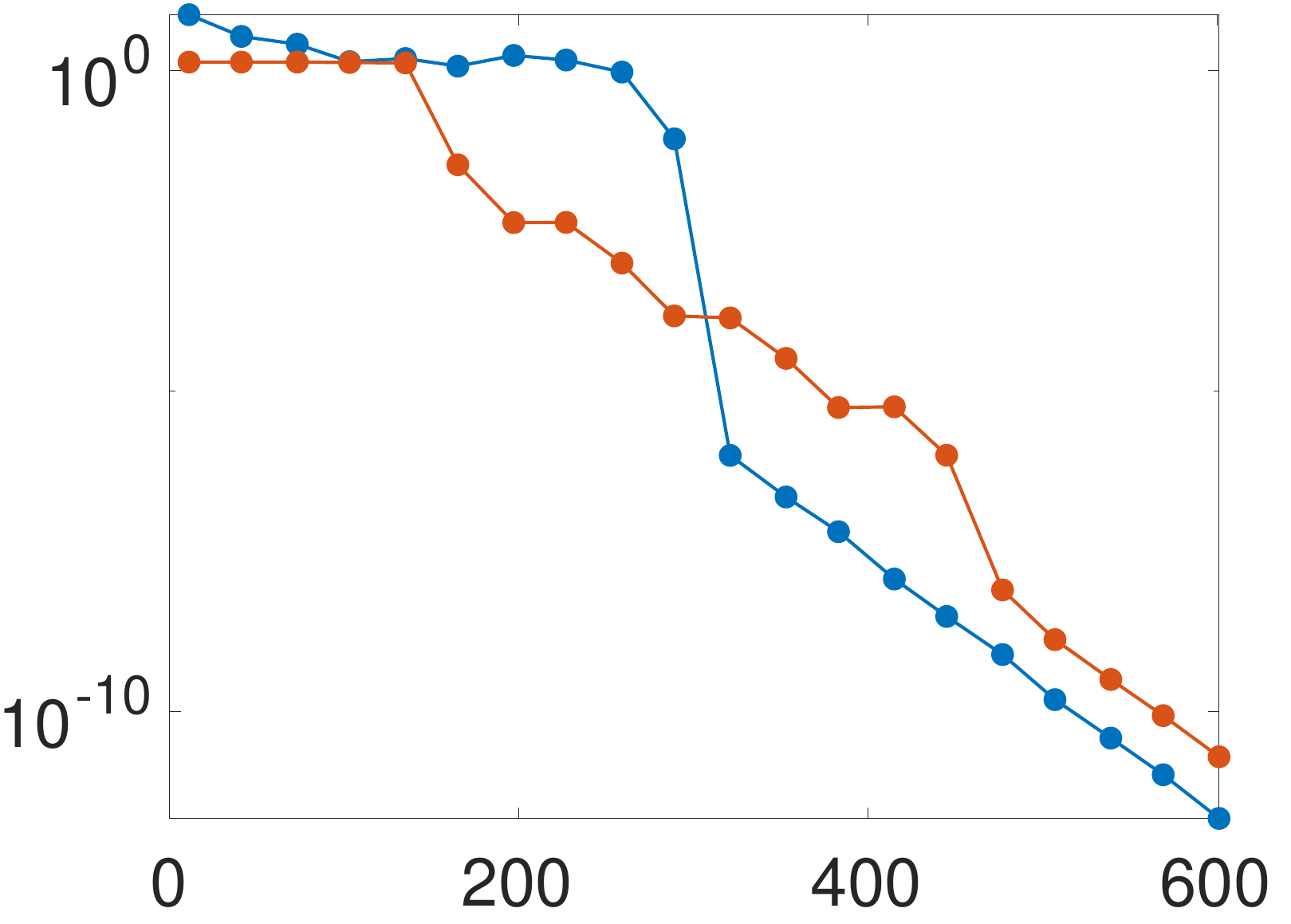}
\put(-135, 40){\rotatebox{90}{{ \small error}}}
\includegraphics[scale = .12]{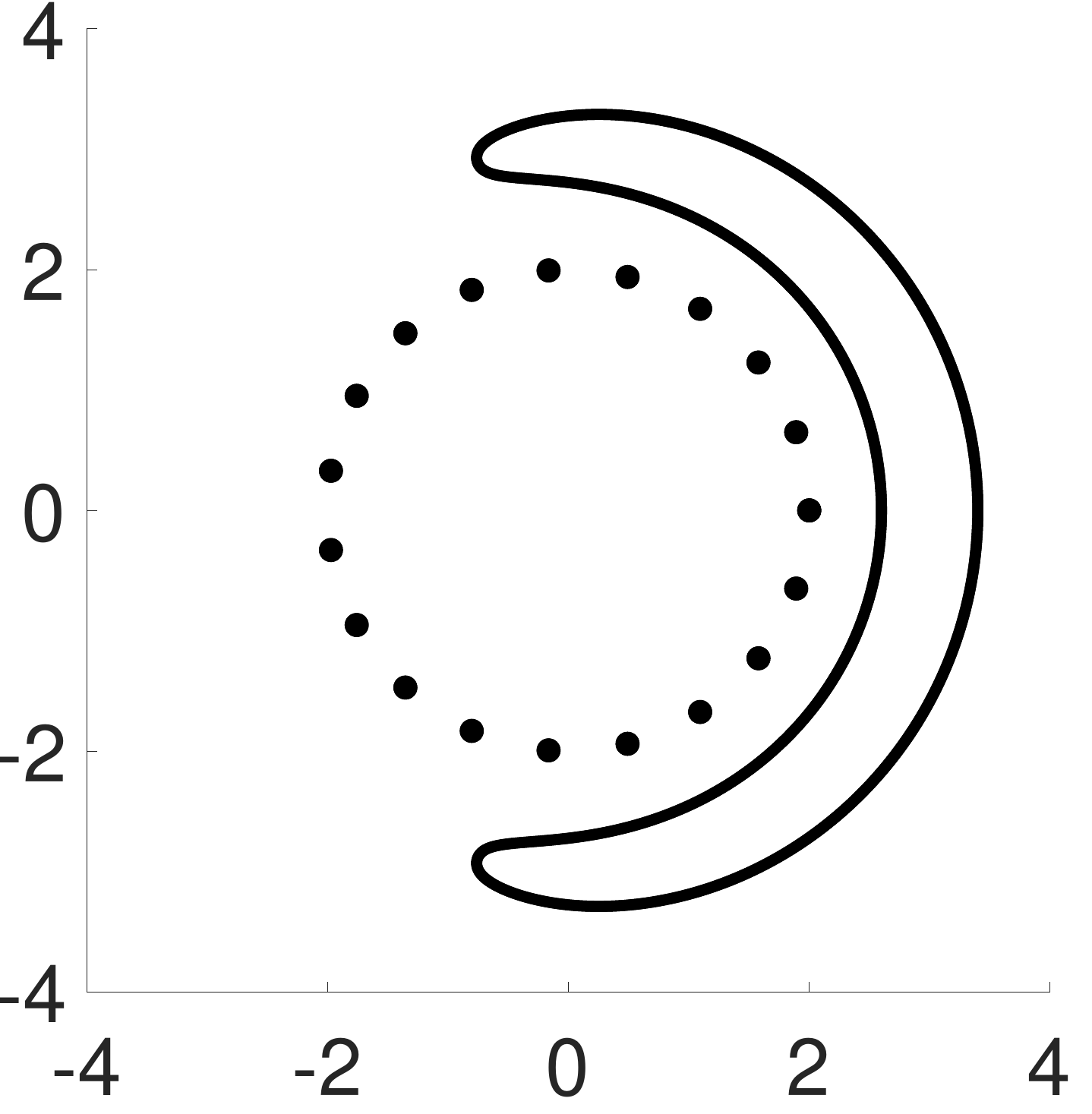}
\newline 
\includegraphics[scale = .155, trim={0cm 0 0cm 0},clip=true]{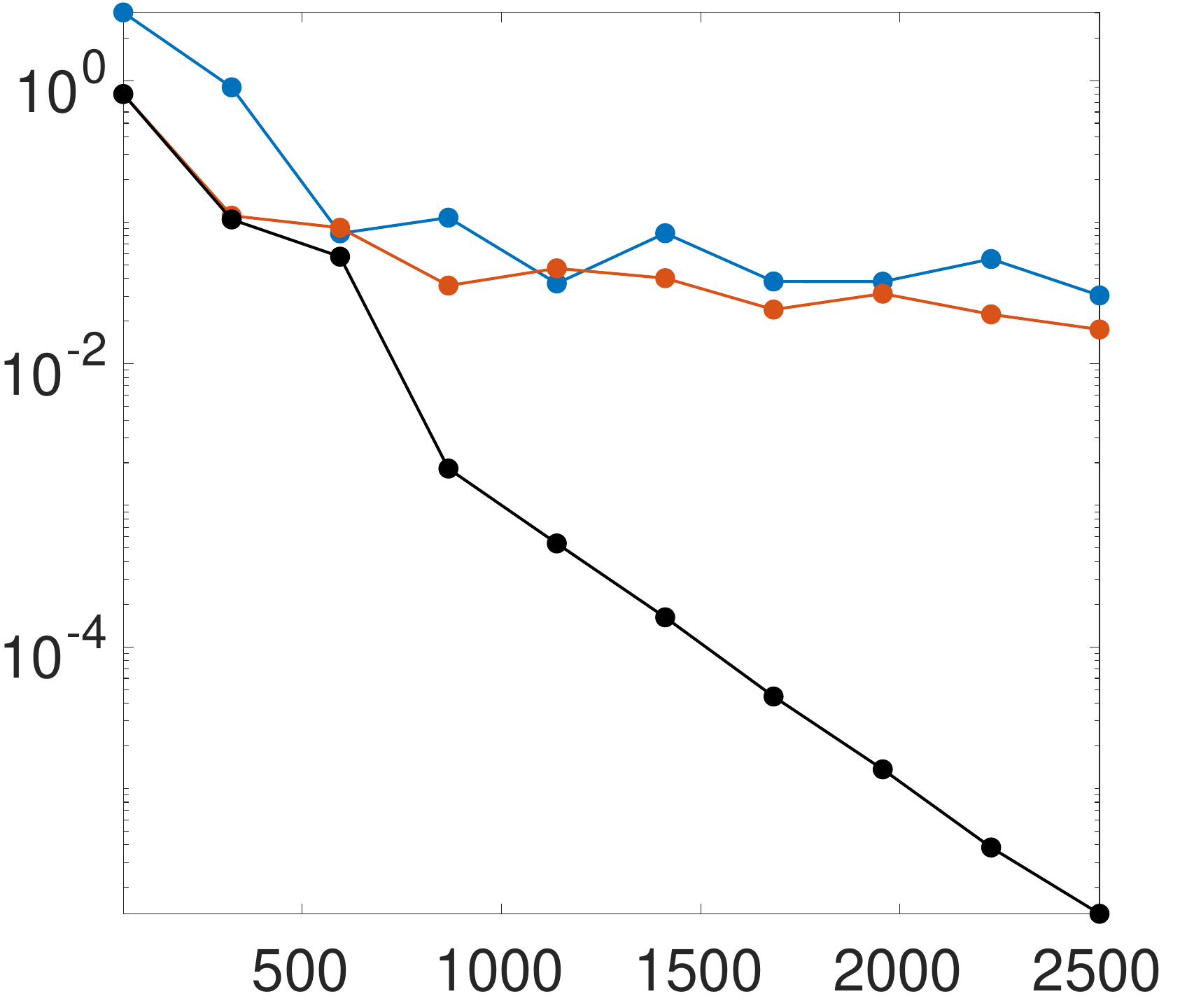}
\put(-135, 50){\rotatebox{90}{{ \small error}}}
\includegraphics[scale = .12, trim={0cm 0 0cm 0},clip=true]{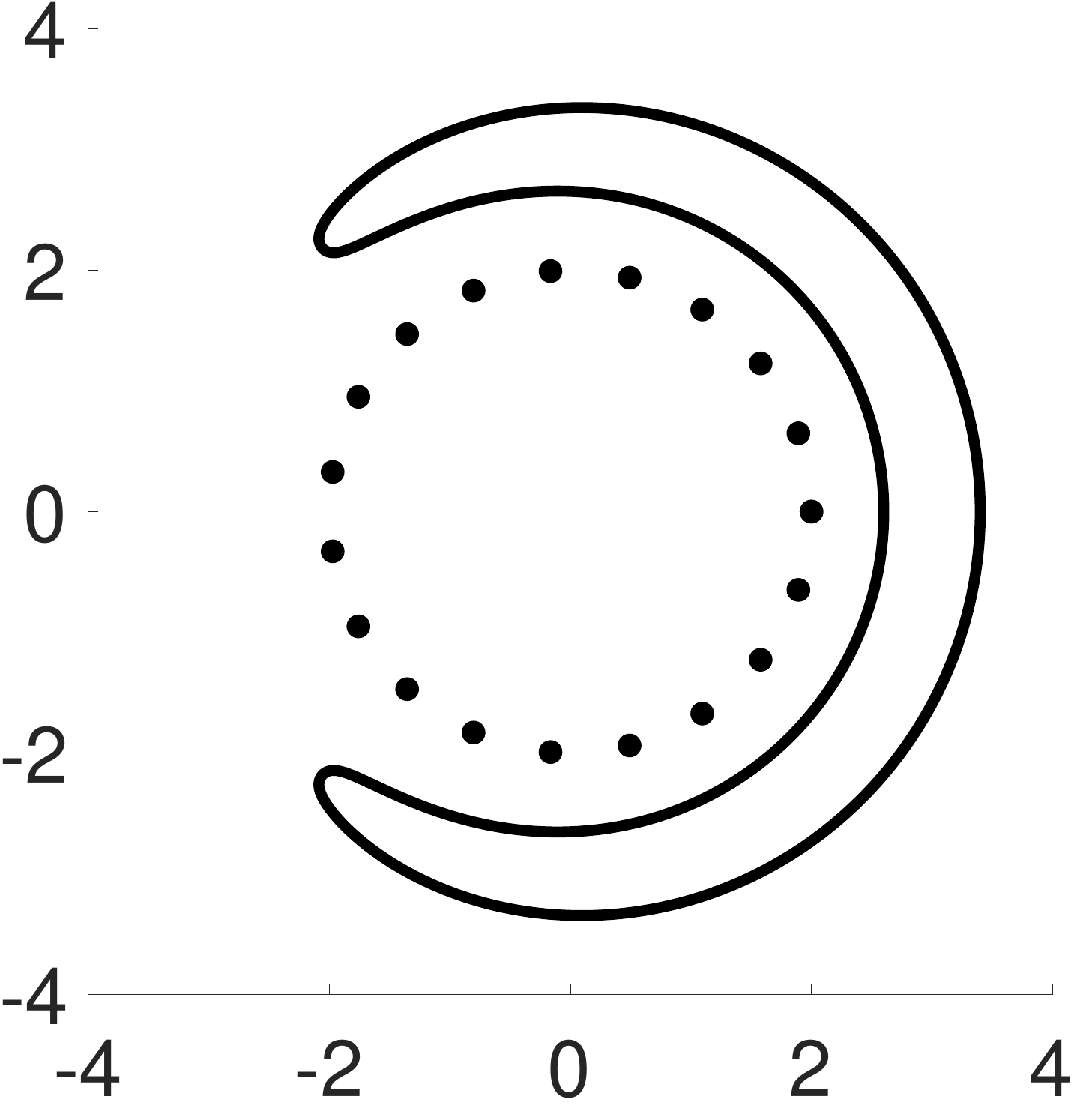}
\caption{ \textit{Top row:} Using the domain on the right, the mean-squared error $\|u_{res}(x,t) -u_m(x,t)\|_2$ over a small set of spatial points $\{x_j = e^{2 \pi i (j-1)/20}\}_{j = 1}^{20}$ (blue dots on right) and time points $\{t_j = 90j\}_{j=0}^{9}$ is plotted on a logarithmic scale against $m$, the number of quadrature points used in computing the inverse Fourier transform, either via Gauss-Legendre quadrature (blue), or via the truncated sinc expansion (red). Here, $u_{res}$ is a highly resolved approximation to the true solution. The convergence of the damping+correction method (black) described in~\Cref{sec:DC} is shown for reference. \textit{Bottom row:} The experiment repeated but with the displayed domain, which induces stronger trapping.}
\label{fig:quadcompare}
 \end{figure}

We can conclude that the regime over which  the sinc-based method is effective and also strongly advantageous over basic quadrature methods is quite limited. This can be seen in~\Cref{fig:quadcompare} (upper row), where the performance of the fast sinc method is compared against that of Gauss-Legendre quadrature in a setting where the trapping is relatively mild. When the trapping effect is marginally stronger, as in the lower row, both methods fail.  In the next section, we describe a way to overcome this limitation via damping so that the sinc-based representation remains useful when there is trapping.
 
\section{A damping+correction scheme for strong trapping} \label{sec:DC}

As indicated in the last section, the decay rate of  $|u(\tilde{x}, t)|$ is an indicator of the convergence rate of $u \to u_m$.  This is closely connected to the distance of poles of the meromorphic extension of $\hat{U}(\tilde{x}, \cdot)$, which come closer to the real line as trapping grows more severe. A natural idea to improve matters is to consider a function with better analyticity properties by perturbing $\omega$ up into the complex plane. An early example of this applied in the context of scattering is found in~\cite{rokhlin1983solution}.  Letting $\omega \in [W_1, W_2]$, the function $\mathcal{F}^{-1}\hat{U}(x, \omega + \delta i)$ corresponds to a solution to the damped wave equation. Using contour integration and the fact that $\hat{U}(x, \cdot)$ is analytic in the upper half of the complex $\omega$-plane, we can use this function to recover the undamped solution without explicitly having to sample $\hat{U}(x, \cdot)$ on the real line.

 We choose $\Gamma = \Gamma_{0} \cup \Gamma_{cR} \cup \Gamma_{\delta} \cup \Gamma_{cL}$ to form a rectangle. Specifically, $\Gamma_0$ is the line segment $[W_1, W_2]$ traversed left to right, $\Gamma_{cR}$ is the line segment $[W_2, W_2+\delta i]$ traversed bottom to top, $\Gamma_\delta$ is the line segment $[W_1 + \delta i, W_2 + \delta i]$ traversed left to right, and $\Gamma_{cL}$ is the line segment $[W_1, W_1 + \delta i]$ traversed top to bottom. Letting $I_a$ denote the contour integral over $\Gamma_a$ we conclude that~\eqref{eq:iFT} is equivalent to  $$I_{0}= I_{\delta} -I_{cR} - I_{cL}.$$  We will refer to $I_{cR}$ and $I_{cL}$ as \textit{correction terms} and $I_{\delta}$ as the \textit{damped term}. 

 \subsection{Evaluating the damped term}
\begin{figure}[h!]
    \centering
      \begin{minipage}{.40\textwidth} 
 \centering
  \begin{overpic}[width=\textwidth]{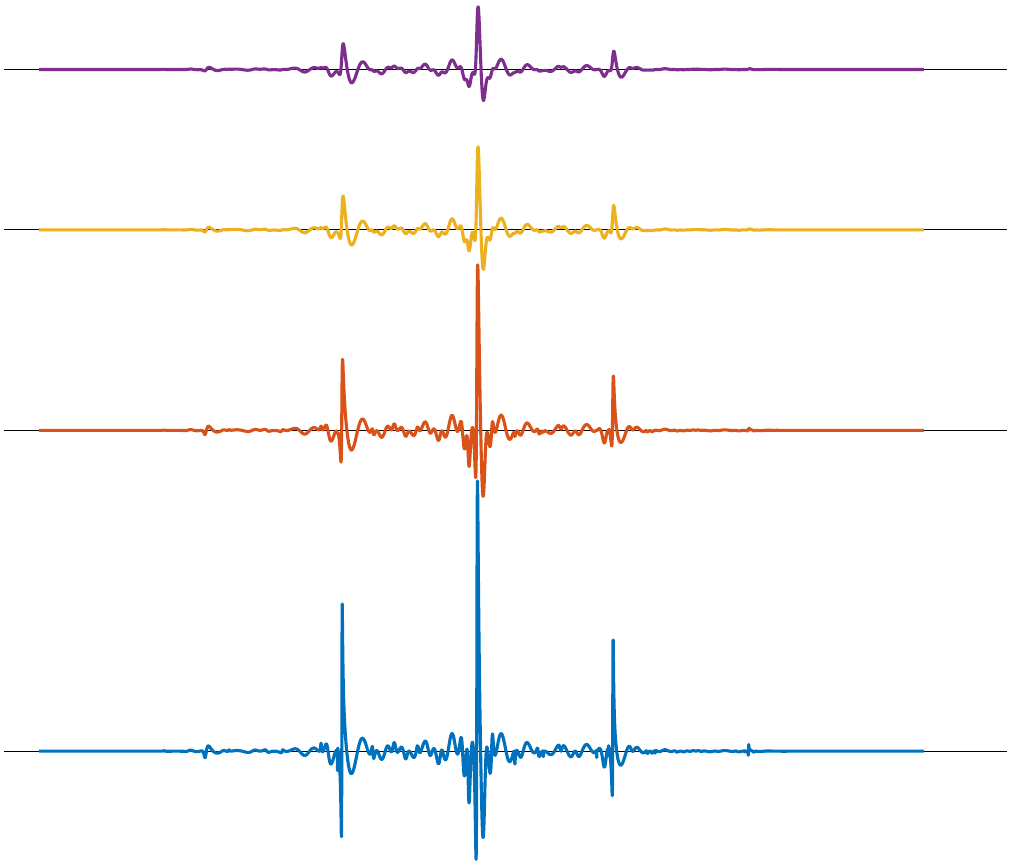}
   \put(10, 14){{\small $\delta = 0$}}
   \put(10, 45){{\small $\delta = .005$}}
   \put(10, 65){{\small $\delta = .015$}}
   \put(10, 80){{\small $\delta =  .02$}}
\end{overpic}
\end{minipage}
	\hspace{.3cm}
\begin{minipage}{.45\textwidth} 
 \centering 
  \begin{overpic}[width=\textwidth]{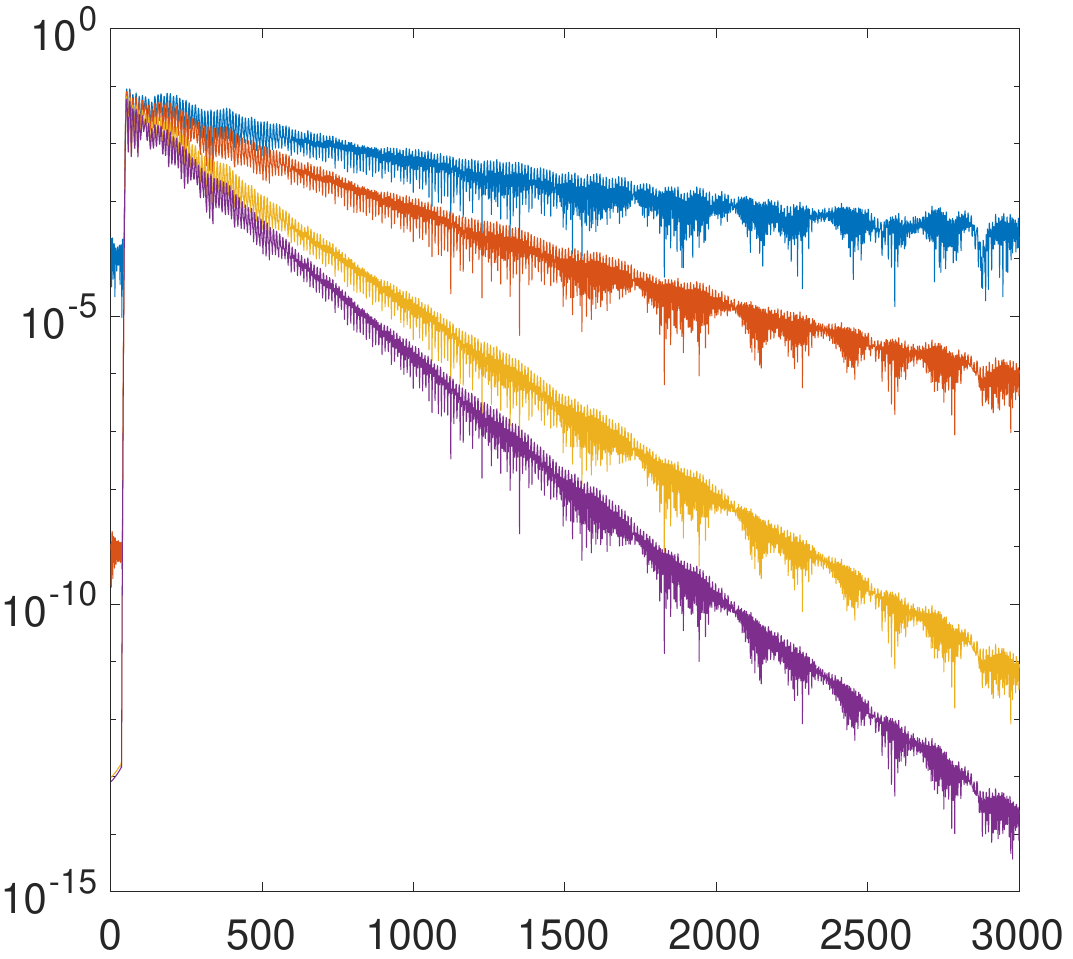}
    \put(50, 75){\rotatebox{-10}{{ \small $\delta = 0$}}}
    \put(61, 64){\rotatebox{-12}{{ \small $\delta = .005$}}}
    \put(61, 49){\rotatebox{-32}{{ \small $\delta = .015$}}}
    \put(58, 29){\rotatebox{-36}{{ \small $\delta = .02$}}}
     \put(-5, 30){\rotatebox{90}{{ \small magnitude}}}
 \put(45, -5){index}
\end{overpic}
\end{minipage}
\caption{ \textit{Left:} Plots of $\hat{U}(\tilde{x}, \omega + \delta i)$ for $\delta = 0$ (blue), $\delta \approx .005 $ (red), $\delta \approx .015 $ (yellow),  and $\delta \approx .02$ (purple) over $\omega \in [ .7, 16]$, with $\tilde{x}$ selected as the center of the $C$-shaped curve in \Cref{fig:poleplot}. The plots are stacked vertically for illustrative purposes. \textit{Right:} The magnitudes of the normalized Fourier coefficients for $\hat{U}(\tilde{x}, \omega + \delta_j i)$ for each $\delta_j$. Faster decay implies that fewer terms are needed to compute the sum in~\eqref{eq:complexsincsum} and thereby evaluate $I_{\delta}$.}
\label{fig:smoothing_slices}
 \end{figure}
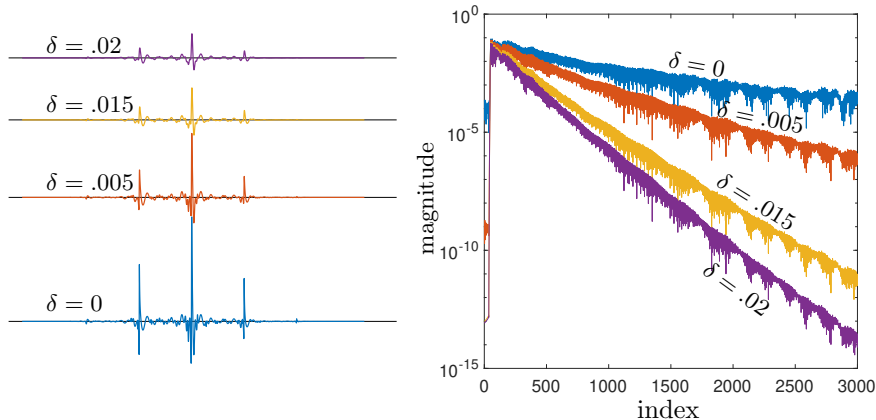

 The damped term can be written as 
 $$ \int_{W_1}^{W_2}  \hat{U}(x, \omega + \delta i) e^{-i (\omega + \delta i) t} {\rm d} \omega, $$
 which simplifies to 
  $$ e^{\delta t} \int_{W_1}^{W_2}  \hat{U}(x, \omega + \delta i) e^{-i\omega t} {\rm d} \omega. $$
  The integral in this expression is of the same form as~\eqref{eq:iFT} and so we can apply the truncated sinc expansion method to evaluate it, which is especially useful because $t$ may be very large so that the integral is highly oscillatory. Moreover, with the truncated sinc expansion, we can use fast transforms to reduce the overall cost of evaluations. This results in an approximation \hbox{$\tilde{u}_m(x,t) \approx \tilde{u}(x, t)$}, where $\tilde{u}(x,t) = \mathcal{F}^{-1} \hat{U}(x, \omega + \delta i)$ and
  \begin{equation} \label{eq:complexsincsum}
   \tilde{u}_m(x,t) =\frac{P}{2 \pi} e^{-i t(P/2 + W_1 + i\delta)} \sum_{j =1}^m (-1)^j \tilde{c}_j(x)   {\rm sinc}\left( \tfrac{Pt}{2 \pi} - j \right).
  \end{equation}
  Here, $\{\tilde{c}_j(x)\}_{j = 1}^m$ are the positively-indexed Fourier coefficients of $\hat{U}(x, \omega + \delta i)$. These coefficients will decay more rapidly than those of $\hat{U}(\tilde{x}, \cdot)$ so that $\tilde{u}_m \to u$ more rapidly than $u_m \to u$. While this improvement may seem modest when $\delta$ is small, it is enough to make a substantial difference in the number of required Helmholtz solves. \Cref{fig:smoothing_slices} illustrates this fact. The left picture shows  profiles of the real part of solutions $\hat{U}(\tilde{x}, \omega + i \delta)$ over $\omega \in [0.7, 16]$ with $\tilde{x} = [0, 0]^T$ selected as the center of the C-shaped keyhole scatterer shown in \Cref{fig:poleplot}. The right image shows the magnitude of the Fourier coefficients of $\hat{U}(\tilde{x}, \omega + \delta i)$ for each $\delta$.

 \subsection{Evaluating the correction terms} The damped terms are not oscillatory integrals and can be evaluated using Gauss-Legendre quadrature. They become increasingly easier to approximate as $\delta$ grows larger, so a natural impulse is to push $\delta$ to be as large as possible. However, exponential growth in the correction terms will cause numerical overflow if $\delta$ is too large, imposing a $T$-dependent limit on how far up $\delta$ can go. A heuristic for estimating this limit is that $\delta$ should not exceed $L(T) = \tau_{mach} \ln(2)/(150T)$, where where $2^{\tau_{mach}}$ is the largest representable floating point number for a given machine. For IEEE 754 double-precision, $\tau_{mach} = 1024$, and $L(T) \approx 6.826 \ln(2)/T$. As an example, for $T = 200$, $L(T) \approx .0237$. For $T = 500$, $L(T) \approx .0095$.

 The heuristic comes from roughly bounding the terms in a $150$ point Gauss-Legendre quadrature estimate of $\int_{0}^\delta e^{\omega T} { \rm d} \omega$. In practice, we adapt the number of quadrature points for evaluating $I_{cL}$ and $I_{cR}$ to the choice of $T$ and $\delta$. In the experiments in \Cref{sec:num}, the number of quadrature points used to evaluate the correction terms ranges between $20$ and $150$.

\subsection{Acoustic scattering for complicated domains}

\begin{algorithm}[t!]
\caption{Damped+correction method for acoustic scattering with trapping}
\label{alg:fullplan}
\begin{algorithmic}[1]
\State If not available, compute $\hat{U}_{inc}(x,\omega)$ for relevant $(x, \omega) \in \partial \Omega \times [W_1, W_2]$.  
\State For $j = 1$ to $m$, use~\eqref{eq:mixedBIE} to compute density functions $\{ \phi_1(x), \cdots, \phi_m(x)\}$ corresponding to solutions $\hat{U}(x, \omega_j + i \delta)$.
\State Evaluate $\hat{U}(x_k, \omega_j + i \delta)$ via~\eqref{eq:U_int} at all values $ (x,\omega) \in \{x_1, \cdots x_M\} \times \{\omega_1, \cdots, \omega_m\}$.
\State Compute Fourier coefficients $\{c_{j k}\}_{j =1, k = 1}^{m,M}$ by applying the FFT to $M$ vectors of the form $\hat{U}(x_k, \bm{\omega})$, where  $\bm{\omega} = \begin{bmatrix} \omega_1, \omega_2, \cdots, \omega_m \end{bmatrix}^T$.
\State For each $(x_k,t_\ell) \in \{x_1, \cdots, x_M\} \times \{ t_1, \cdots, t_N\}$, evaluate the integral $I_\delta$ using~\eqref{eq:complexsincsum}. 
\State Evaluate the correction terms $I_{cL}$ and $I_{cR}$ at all $(x_j,t_\ell)$.
\State Set $u(x_k,t_\ell) := I_{\delta}(x_k,t_\ell) -I_{cR}(x_k,t_\ell) - I_{cL}(x_k, t_\ell)$.
\end{algorithmic}
\end{algorithm}

Combining the contour deformation with the fast sinc transform method for handling large $t$ leads to a relatively straightforward and effective scheme for solving the acoustic scattering problem on a variety of challenging domains. We describe in pseudocode in \Cref{alg:fullplan} the workflow for solving~\eqref{eq:IVBP}--\eqref{eq:IBVP3}.

\begin{figure}[t!]
 \centering
 \hspace{-.4cm} \includegraphics[scale = .18, trim={2cm 2cm 3cm 0cm},clip]{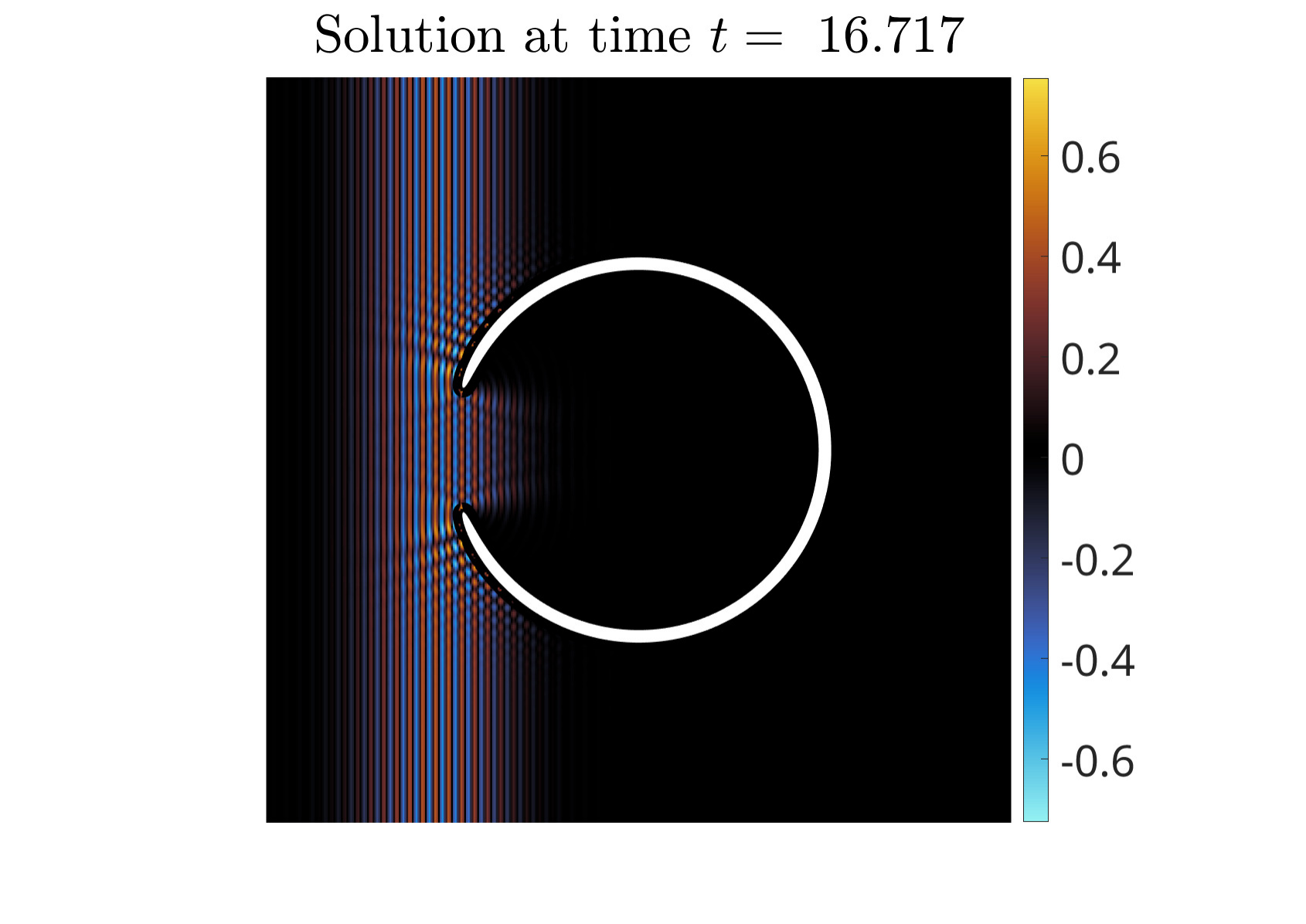}
\includegraphics[scale = .18, trim={2cm 2cm 3cm 0cm},clip=true]{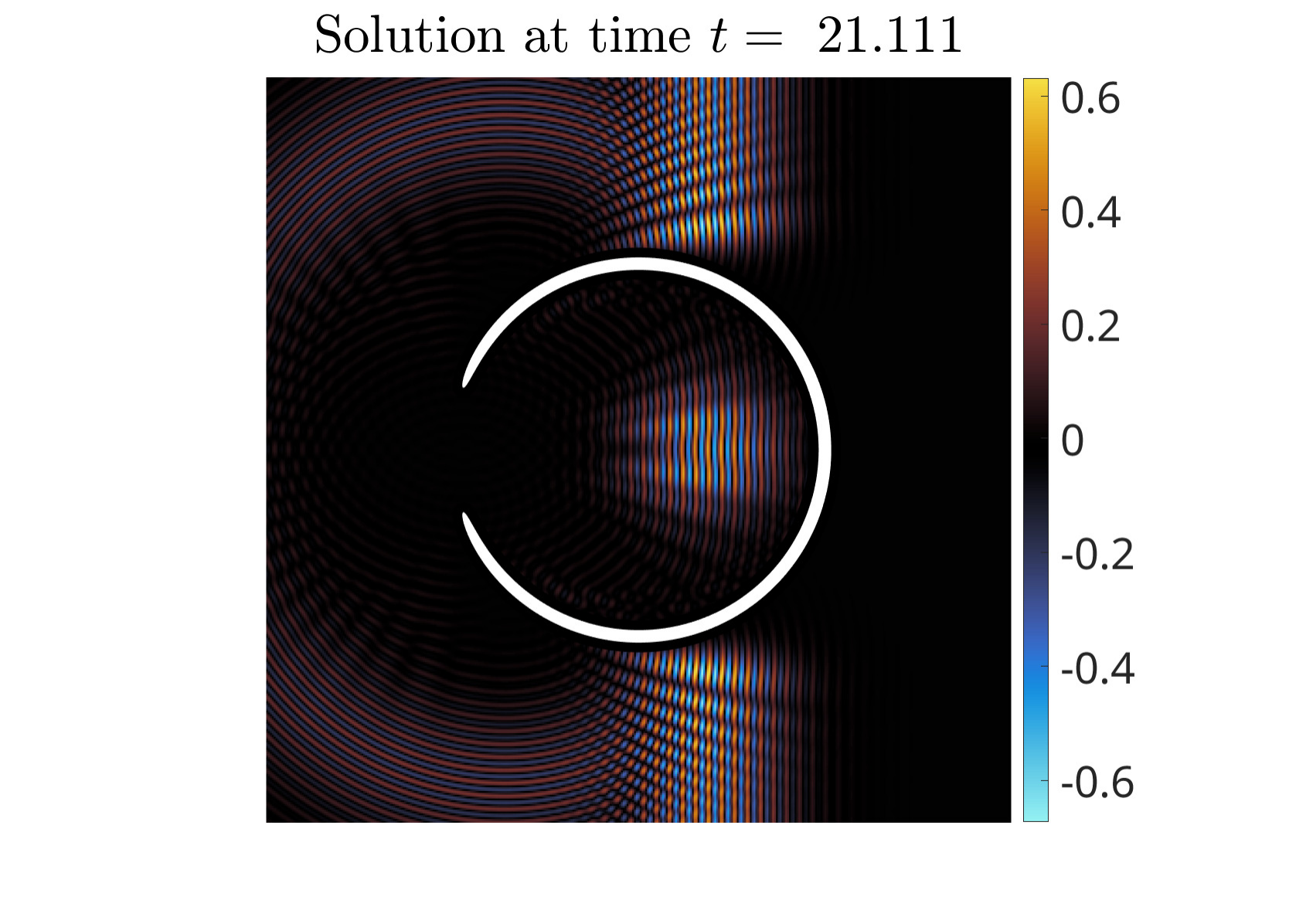}
\includegraphics[scale = .18, trim={2cm 2cm 3cm 0cm},clip=true]{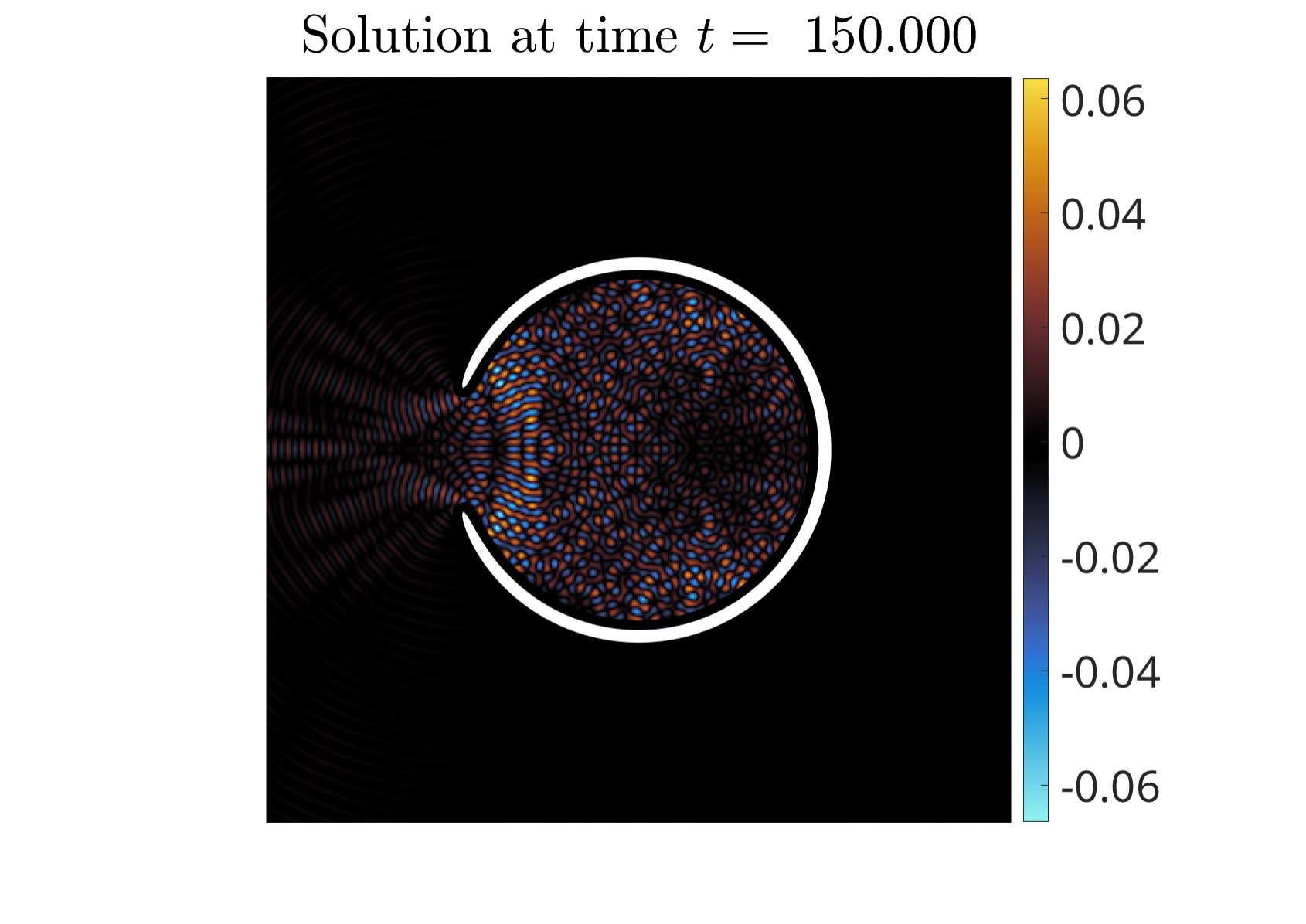}
\caption{The real part of the total field for a solution to the acoustic wave equation is plotted at $t = 16.717$ (left), $t = 21.111$ (center), and $t = 150$ (right).}
  \label{fig:c_curveplots}
\end{figure}

\section{Numerical Results} \label{sec:num}
  In this section we include several experiments that illustrate the capabilities of the solver. Videos of simulations can be found at the link listed in~\cite{HOMEyoutube}.
  We have implemented the solver in MATLAB, with all code publicly available in~\cite{TFAW}. We make use of several existing packages, including Flatiron's FINUFFT library~\cite{barnett2018non} and their implementation of the fast-sinc transform~\cite{lawrenceFS}, the fmmlib2D~\cite{greengard2020fmmlib2d} suite, and the Zetatrap package~\cite{wu2021zeta}. Specialized colormaps used in our visualizations are from the slanCM package~\cite{Slanpack}. Finally, we make use of the so-called \texttt{chop} function developed in Chebfun~\cite{aurentz2017chopping} for checking the resolution quality of~\eqref{eq:sincsum} and enabling the adaptive selection of $m$. To improve the accessibility of our experiments, we include MATLAB code in \Cref{apx:scatter} that generates all of the scatterers we use. 

    \begin{figure}[ht!]
    \centering
    \begin{minipage}{.37\textwidth} 
 \centering
  \begin{overpic}[width=\textwidth]{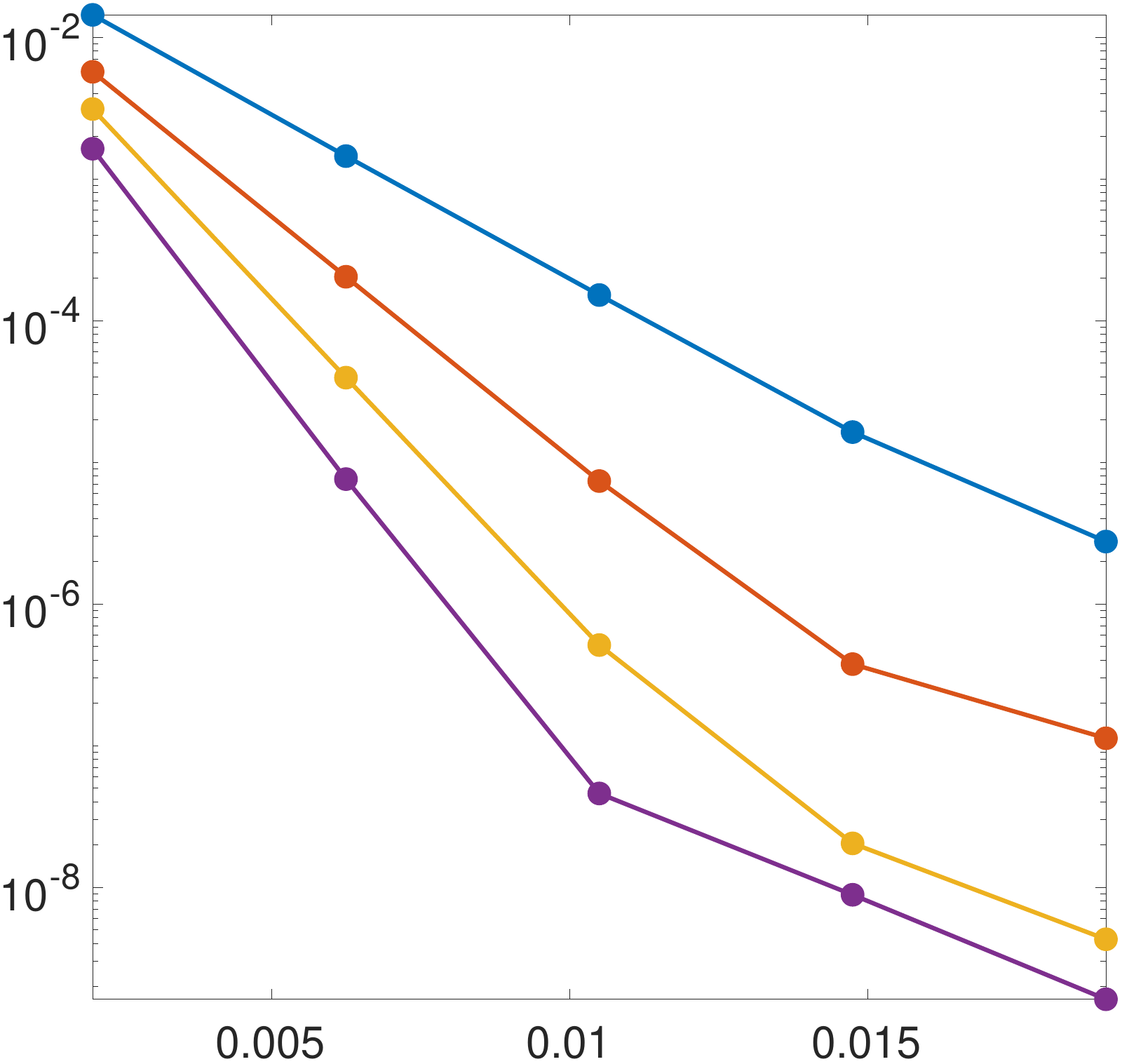}
  \put(55, -4){{\small $\delta$}}
     \put(43, 75){\rotatebox{-25.65}{{ \small $m = 1475$}}}
     \put(40, 65){\rotatebox{-36}{{ \small $m = 2150$}}}
     \put(40, 53){\rotatebox{-43}{{ \small $m = 2825$}}}
     \put(34, 51){\rotatebox{-51}{{ \small $m = 3500$}}}
\end{overpic}
\end{minipage}
	\hspace{.3cm}
\begin{minipage}{.35\textwidth} 
 \centering 
  \begin{overpic}[width=\textwidth]{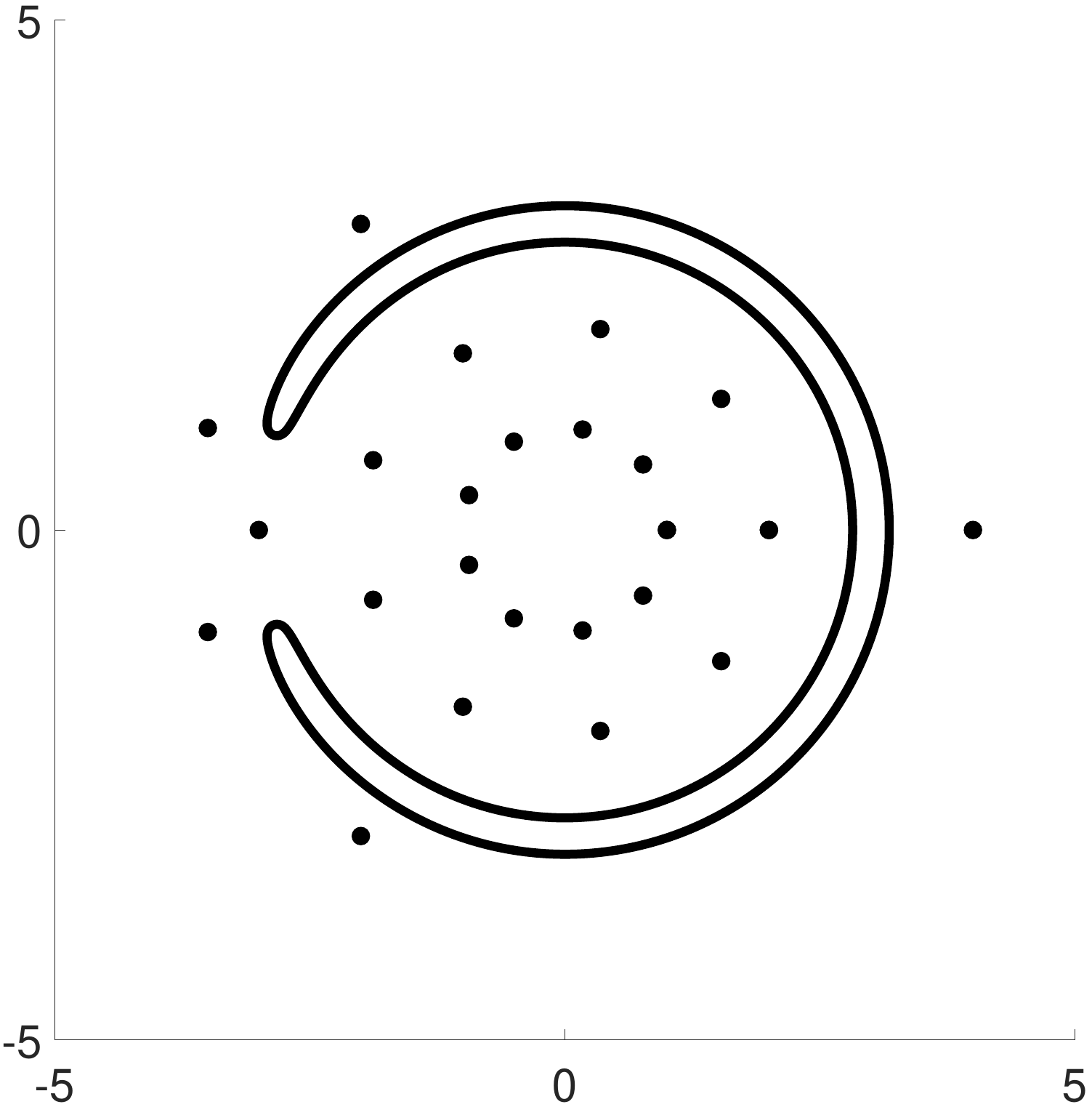}
\end{overpic}
\end{minipage}
\caption{ The $2$-norm error of solutions (relative to a highly resolved solution) computed over three points in time and the 20 spatial points shown on the right (black dots) is plotted on a logarithmic scale against the parameter $\delta$. Each line shows the error behavior for a fixed choice of $m$ as $\delta$ grows. }  \label{fig:c_curve_test}
 \end{figure}
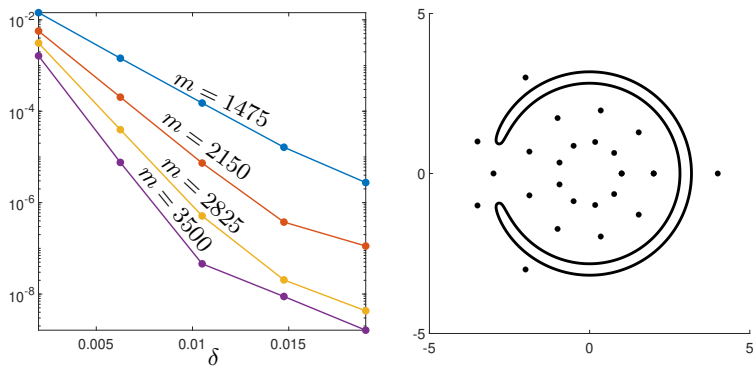

 \subsection{Convergence behavior in a smooth domain} In this example, we examine the behavior of the solution as $m$ and $\delta$ are varied. The problem we consider involves a $C$-shaped scatterer with a smooth boundary. The incident wave is set as in~\eqref{eq:enveloped_osc_inc_wave} with $\omega_0 = 30$, $z_0 = [1,0]^T$, and $t_0 = 20$. The frequency band over which $\hat{U}(x, \cdot)$ is relevant is $[21, 38]$. \Cref{fig:c_curveplots} shows the geometry and the real part of the total field at various points in time. Note that at $t = 150$, the scattered field is still on the order of $10^{-2}$. 

  The images in \Cref{fig:c_curveplots} were produced with a highly resolved approximation to the solution with $m = 4000$ and $\delta = .02$. In \Cref{fig:c_curve_test}, we look at the  error between this solution and solutions computed with different choices of $\delta$ and $m$. We measure the 2-norm error calculated over $20$ points in the domain (see \Cref{fig:c_curve_test}) and five equally-spaced points in time taken over the interval $[5, 150]$.

 \subsection{Multiply-connected domains} 
\begin{figure}[t!]
\centering
\includegraphics[scale = .16, trim={3cm 1cm 1cm 0cm},clip=true]{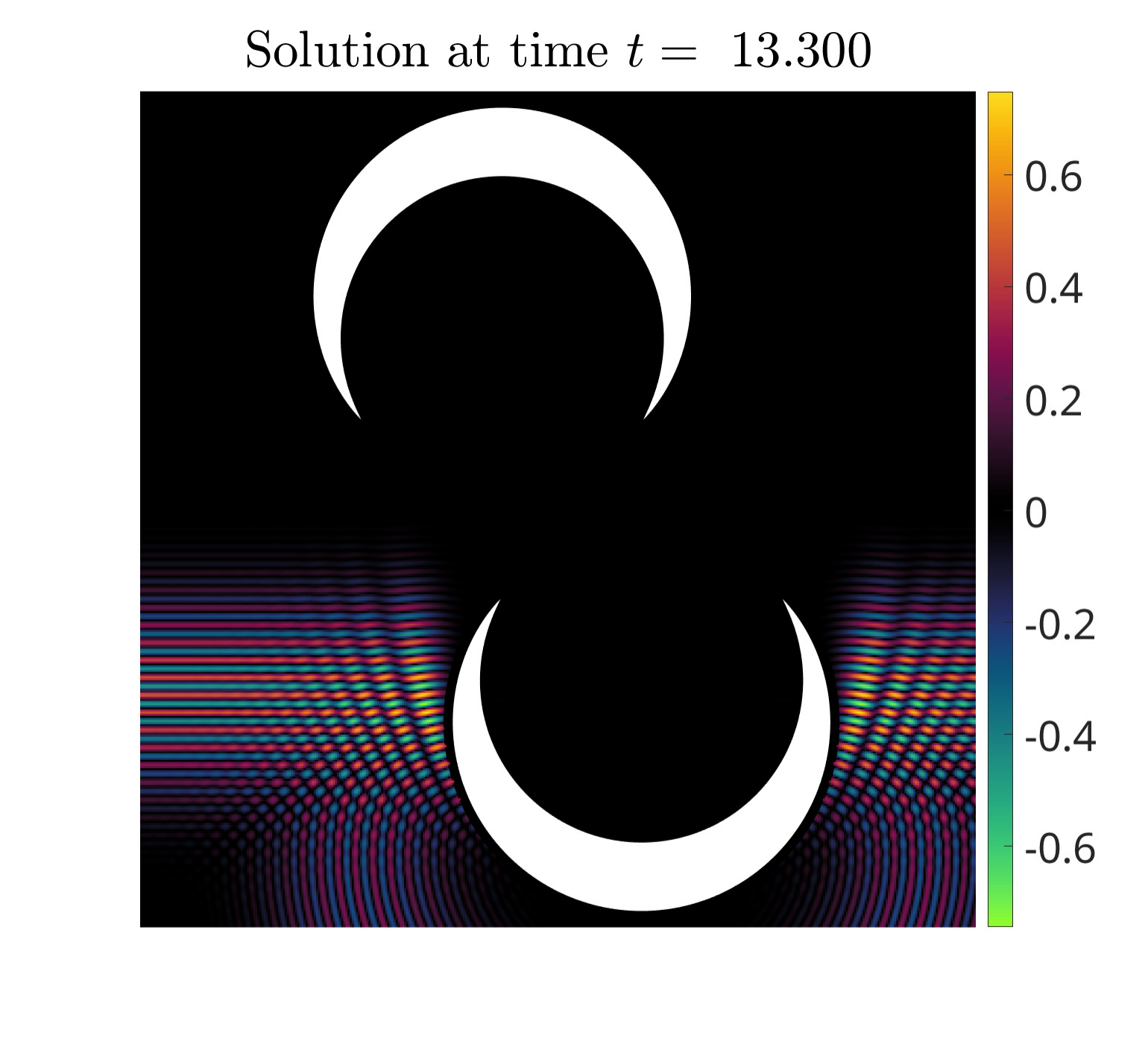}
\includegraphics[scale = .16, trim={3cm 1cm 1cm 0cm},clip=true]{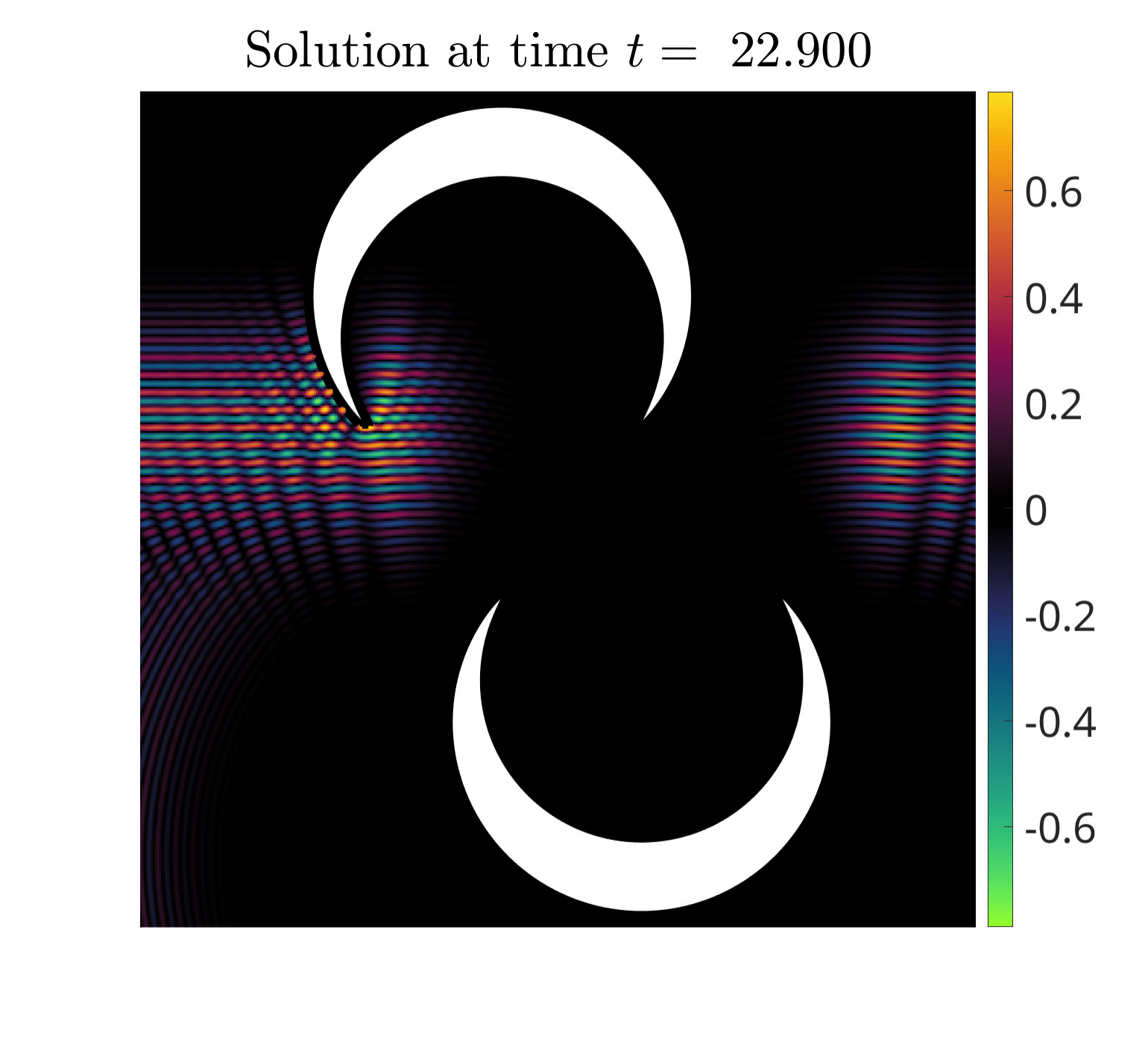}
\includegraphics[scale = .165, trim={3cm 1cm 1cm 0cm},clip=true]{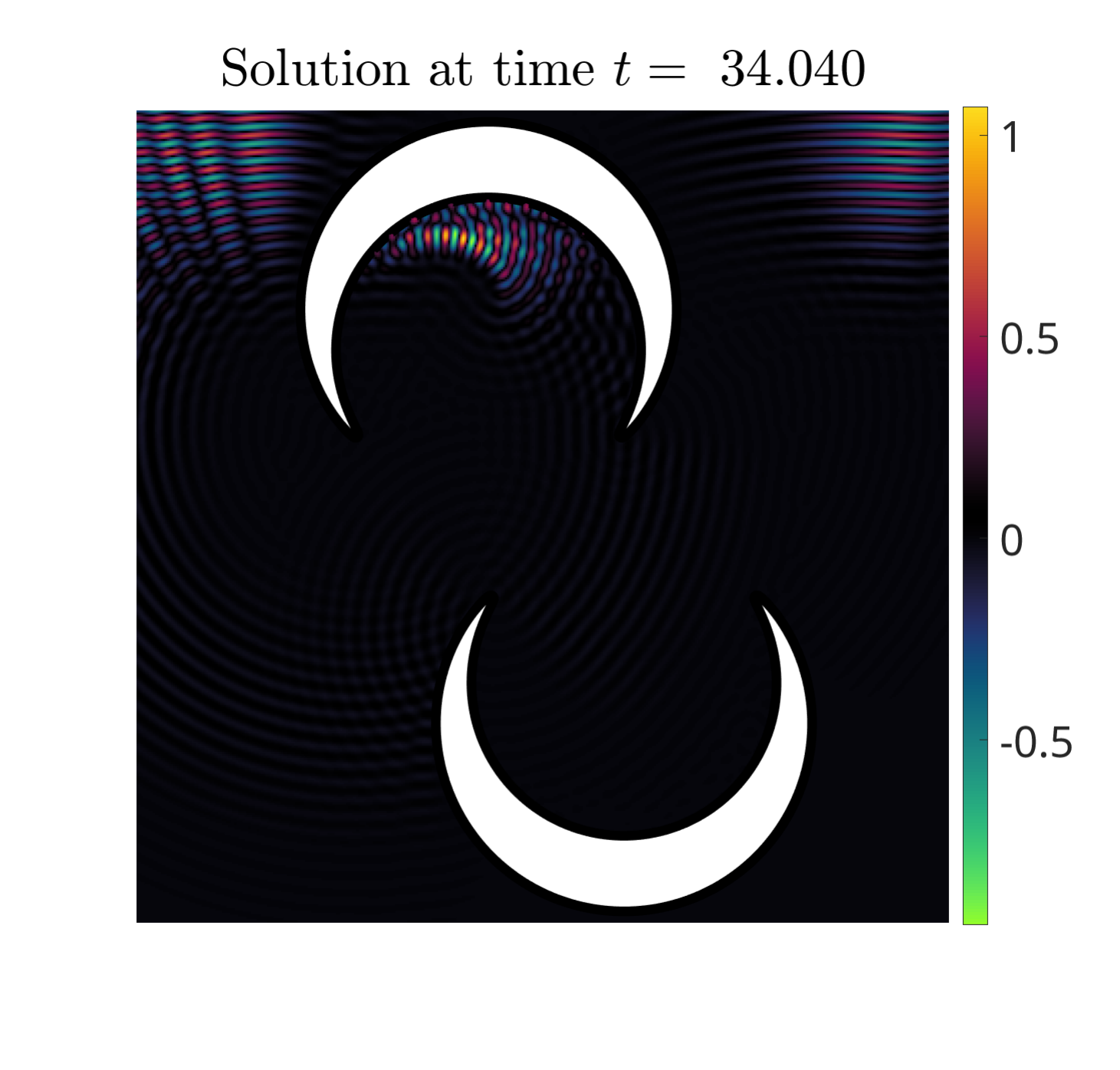}
\includegraphics[scale = .16, trim={3cm 1cm 1cm 0cm},clip=true]{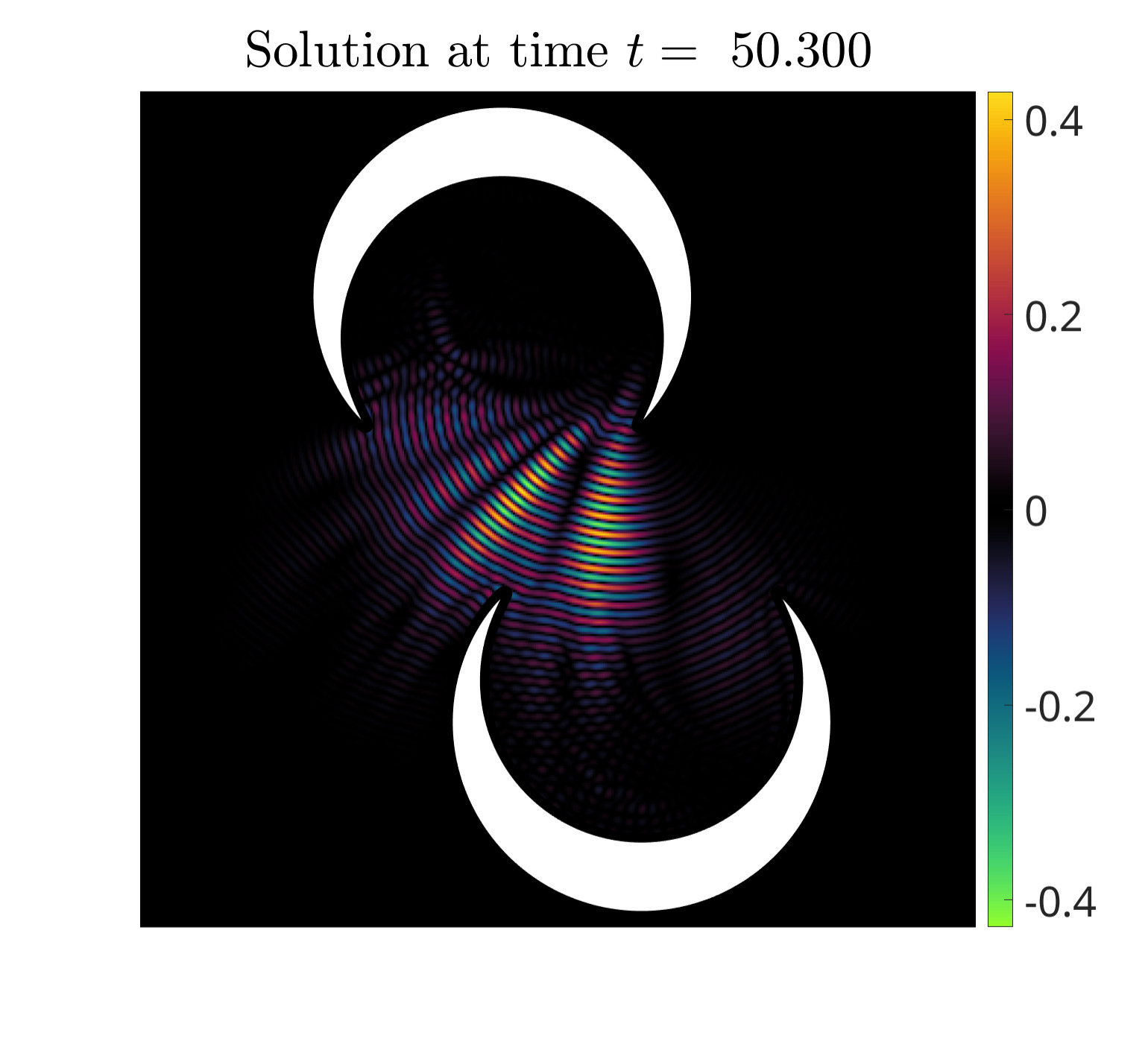}
\includegraphics[scale = .16, trim={3cm 1cm 1cm 0cm},clip=true]{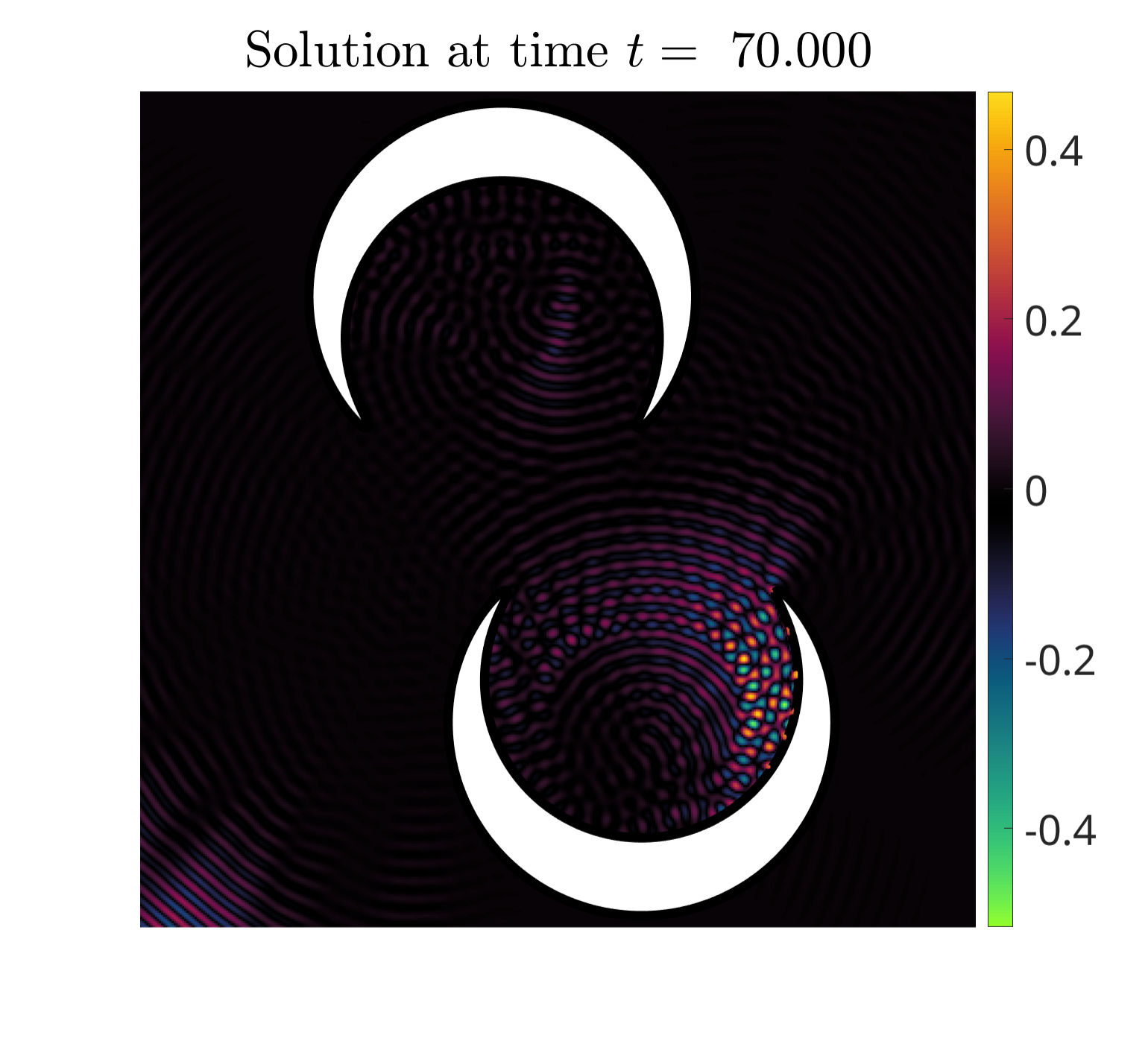}
\includegraphics[scale = .16, trim={3cm 1cm 1cm 0cm},clip=true]{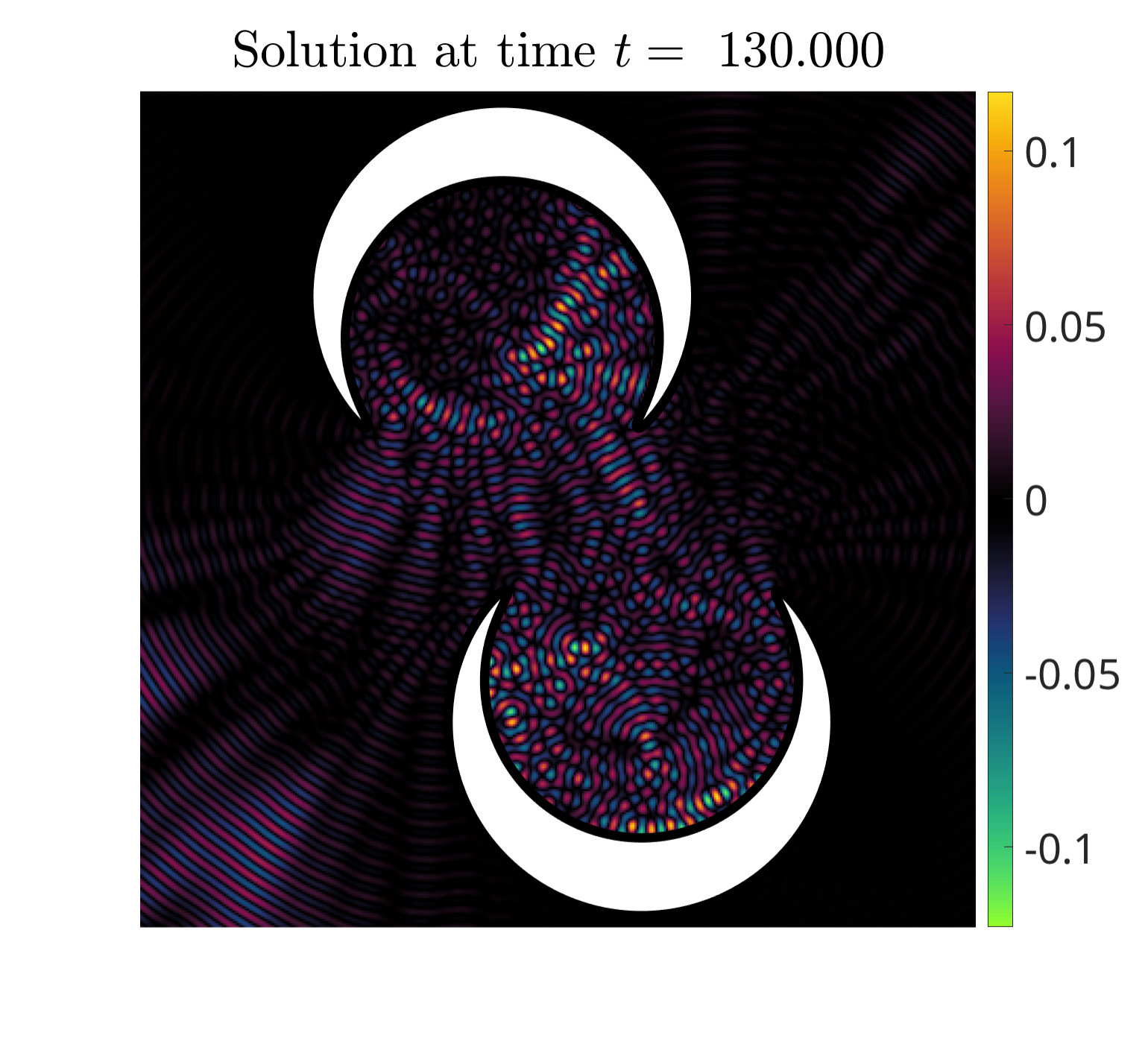}
\caption{The real part of a solution where the scattered wavefield is passed back and forth between two offset crescents. Times included are $t = 13.3, 22.9, 34.04$ (top row, left to right), and $t = 50.3, 70.0, 130.0$ (bottom row, left to right)}
\label{fig:crescent}
 \end{figure}

No special modifications are required to handle multiply-connected domains. In \Cref{fig:crescent} we show a solution involving two offset crescents that ``pass" the wavefield between one-another. For this example the incident wave is chosen as in~\eqref{eq:enveloped_osc_inc_wave} with $\omega_0 = 10$, $z_0 = [0,1]^T$, and $t_0 = 20$. The frequencies over which the problem is relevant are $\omega \in [7, 12]$, and we set $\delta = .02$ with $m = 850$.

  \subsection{Domains with corners} \label{sec:cornerdoms}

  \begin{figure}[h!]
\centering
\includegraphics[scale = .24,trim={2.5cm 2cm 1cm 2cm},clip=true]{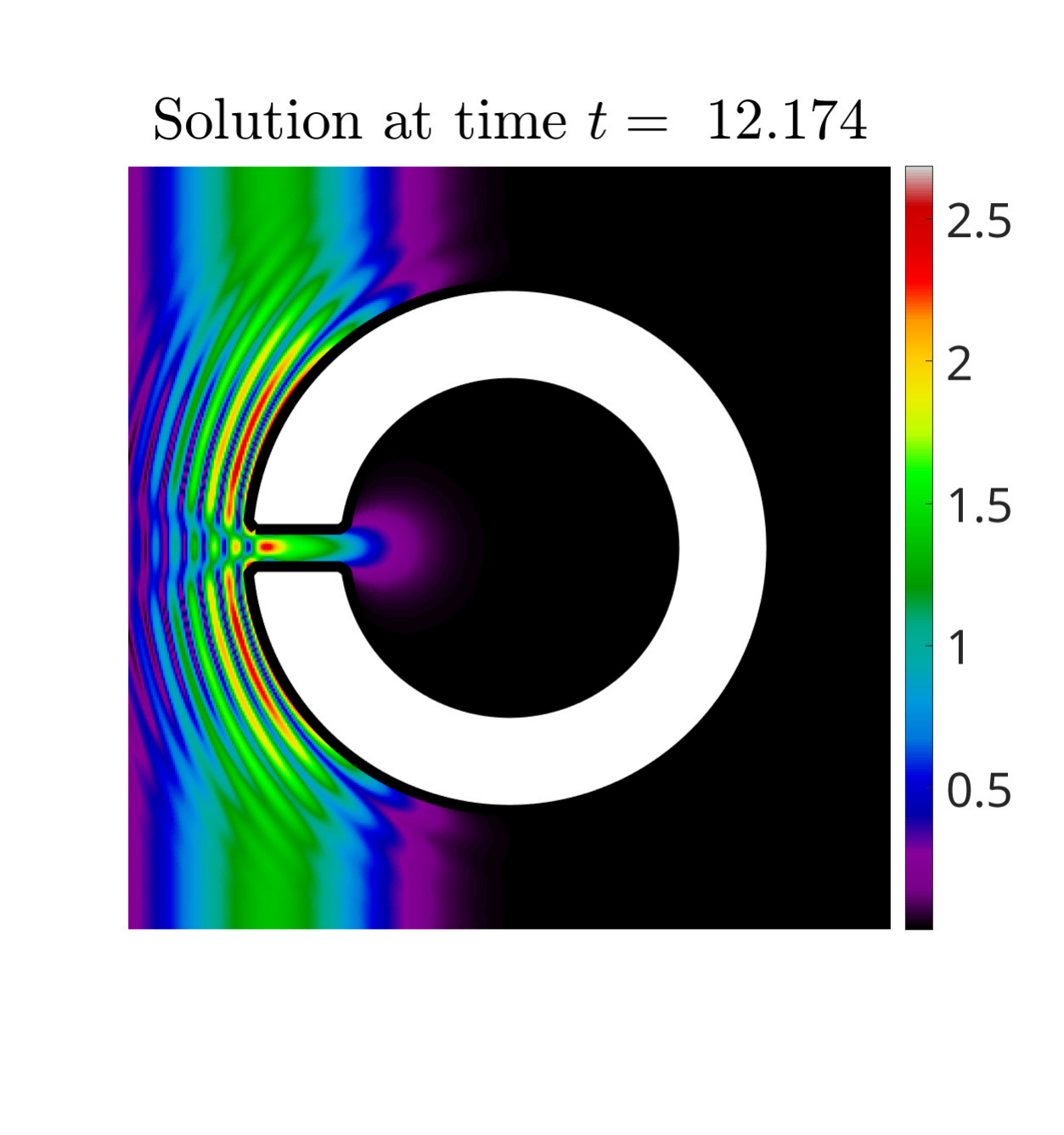}
\includegraphics[scale = .24, trim={2.5cm 2cm 1cm 2cm},clip=true]{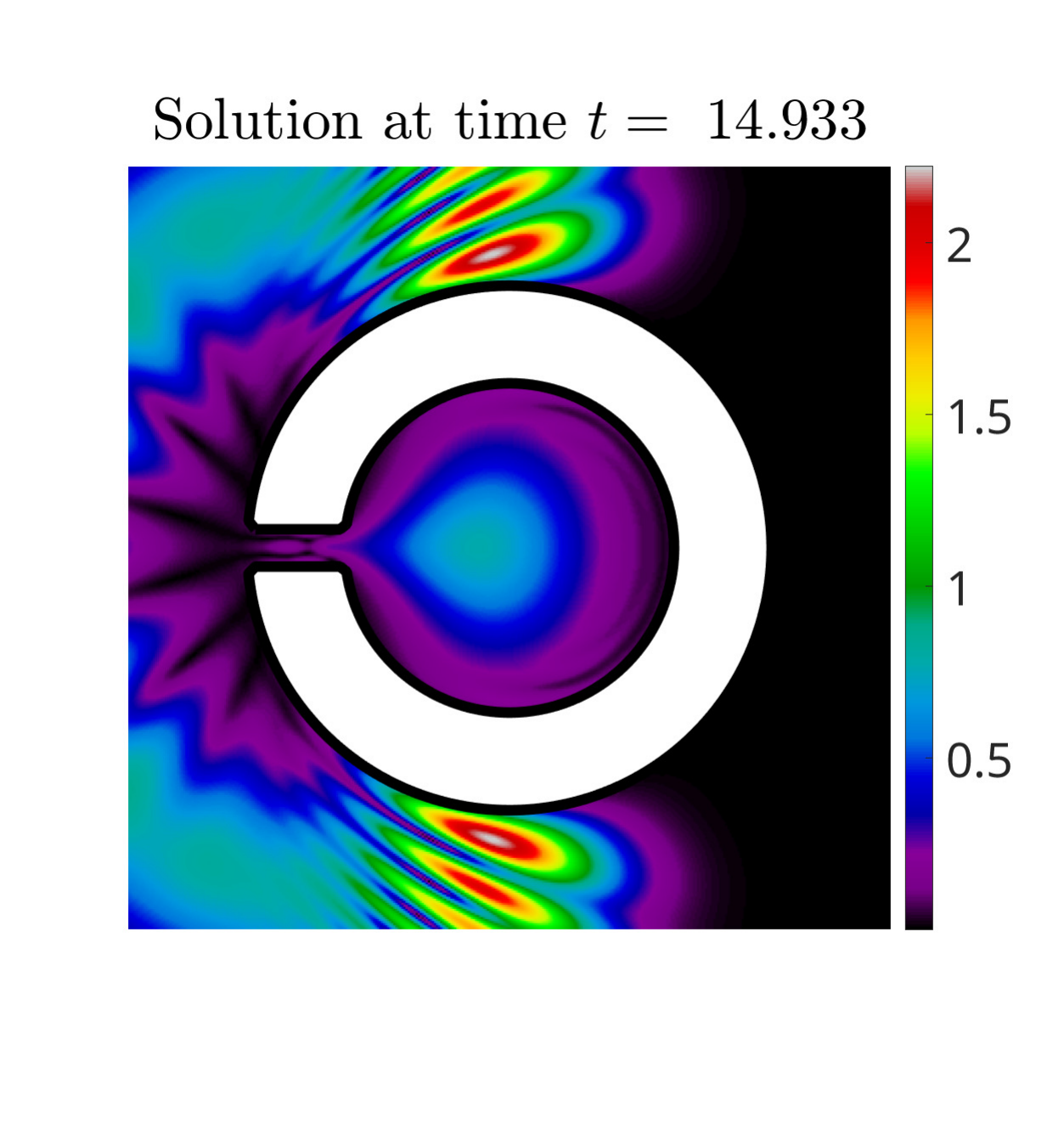}
\includegraphics[scale = .24, trim={2.5cm 2cm 1cm 2cm},clip=true]{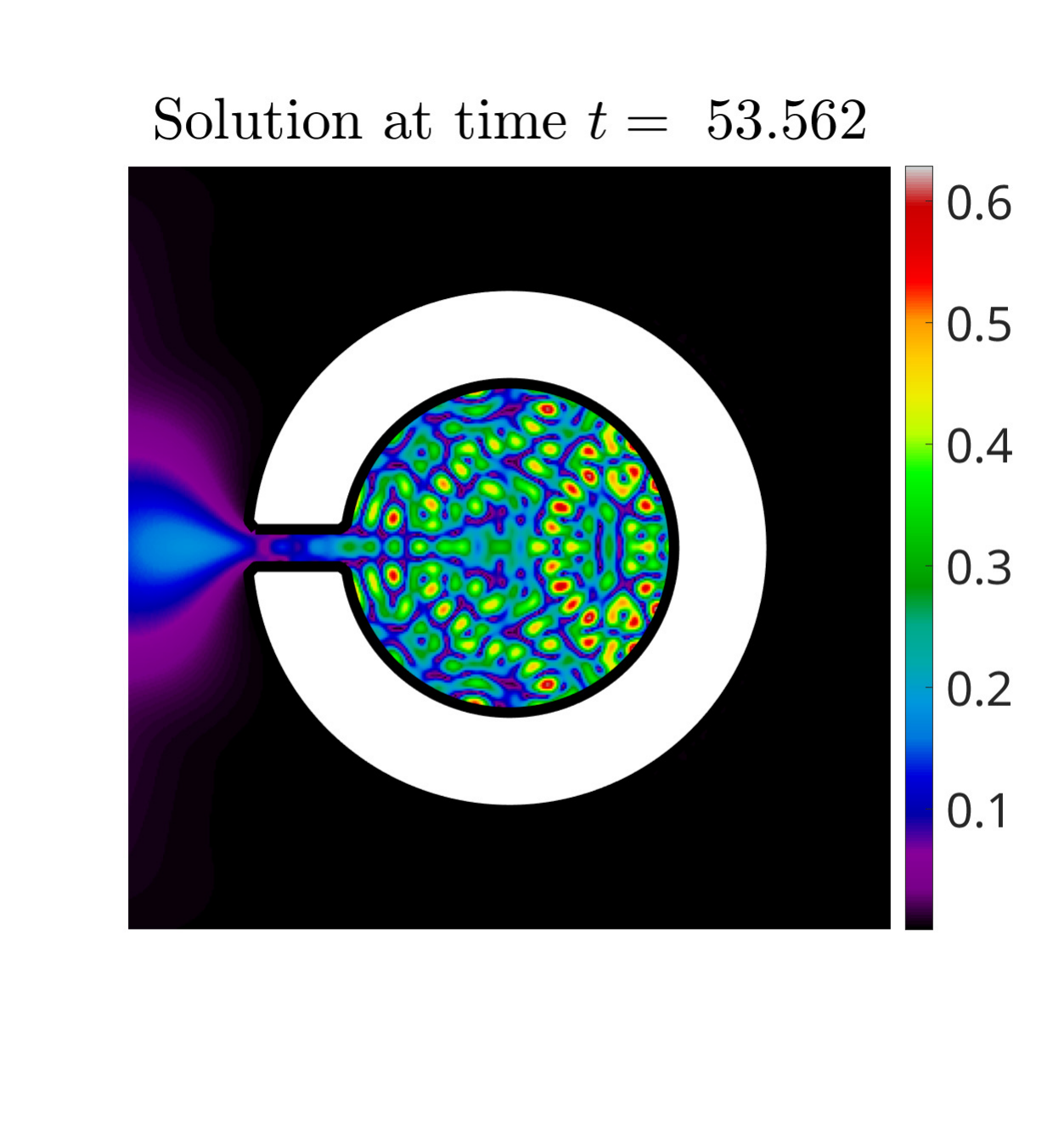}
\vspace{-1cm}
\caption{ Magnitude of total field in a keyhole domain at times $t = 12.174$ (left), $t = 14.933$ (center), $t = 53.562$ (right).}
\label{fig:cornerskey}
\end{figure}

Time-frequency-based methods possess two major potential advantages over methods that discretize directly in time. First, the former is not burdened by CFL-like restrictions that can provide highly restrictive bounds on the maximal time-step allowed, in particular in situations that involve mesh refinement. Second, the spatial discretization for solving a 2D problem only involves resolution along a 1D curve (the boundary of $\Omega$). The damping+correction fast-sinc method is effective for handling domains with both corners and trapping, which is a regime in which traditional time-stepping methods seriously struggle.  

\begin{figure}[h!]
\centering
\includegraphics[scale = .22, trim={0cm 0 0cm 0},clip=true]{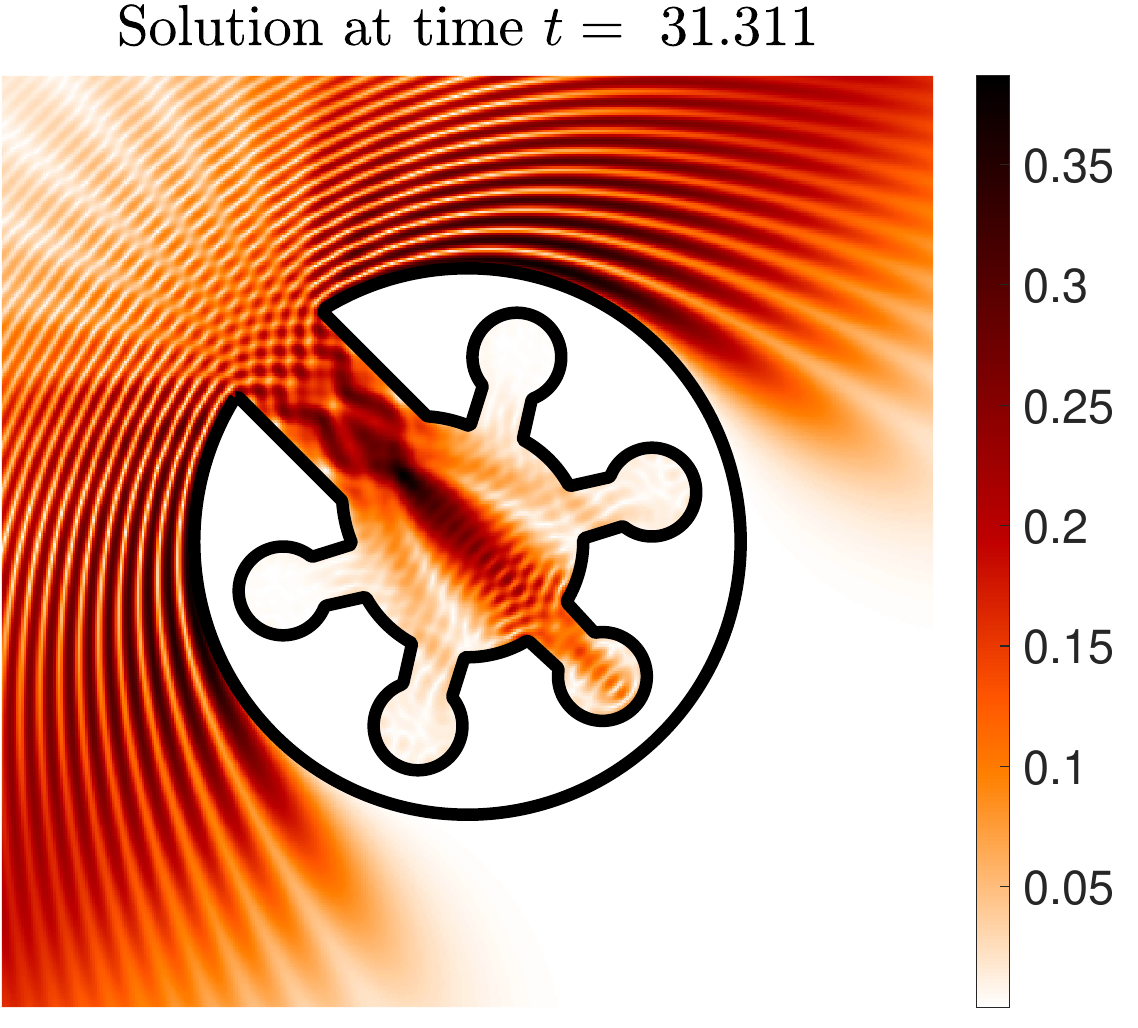}
\includegraphics[scale = .22, trim={0cm 0 0cm 0},clip=true]{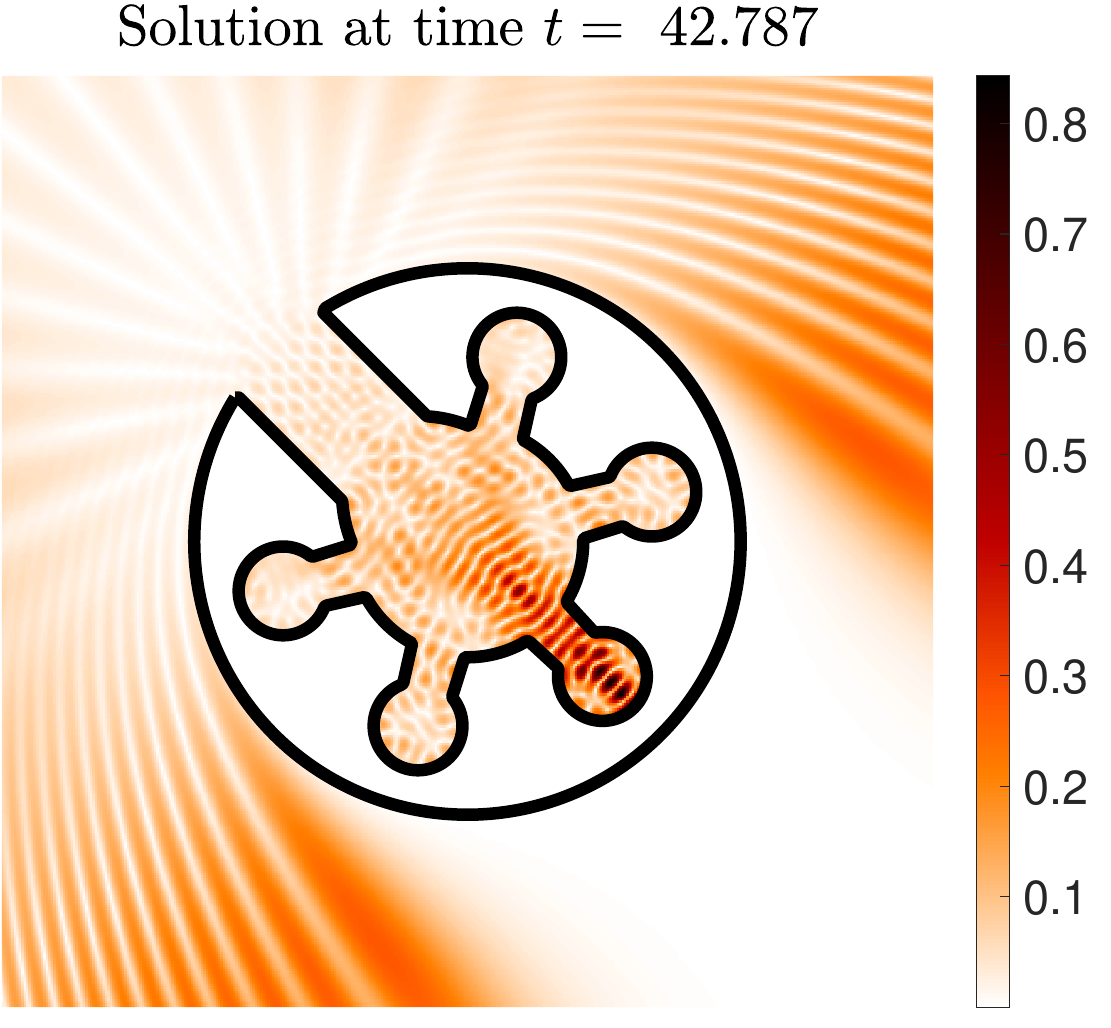}
\includegraphics[scale = .22, trim={0cm 0 0cm 0},clip=true]{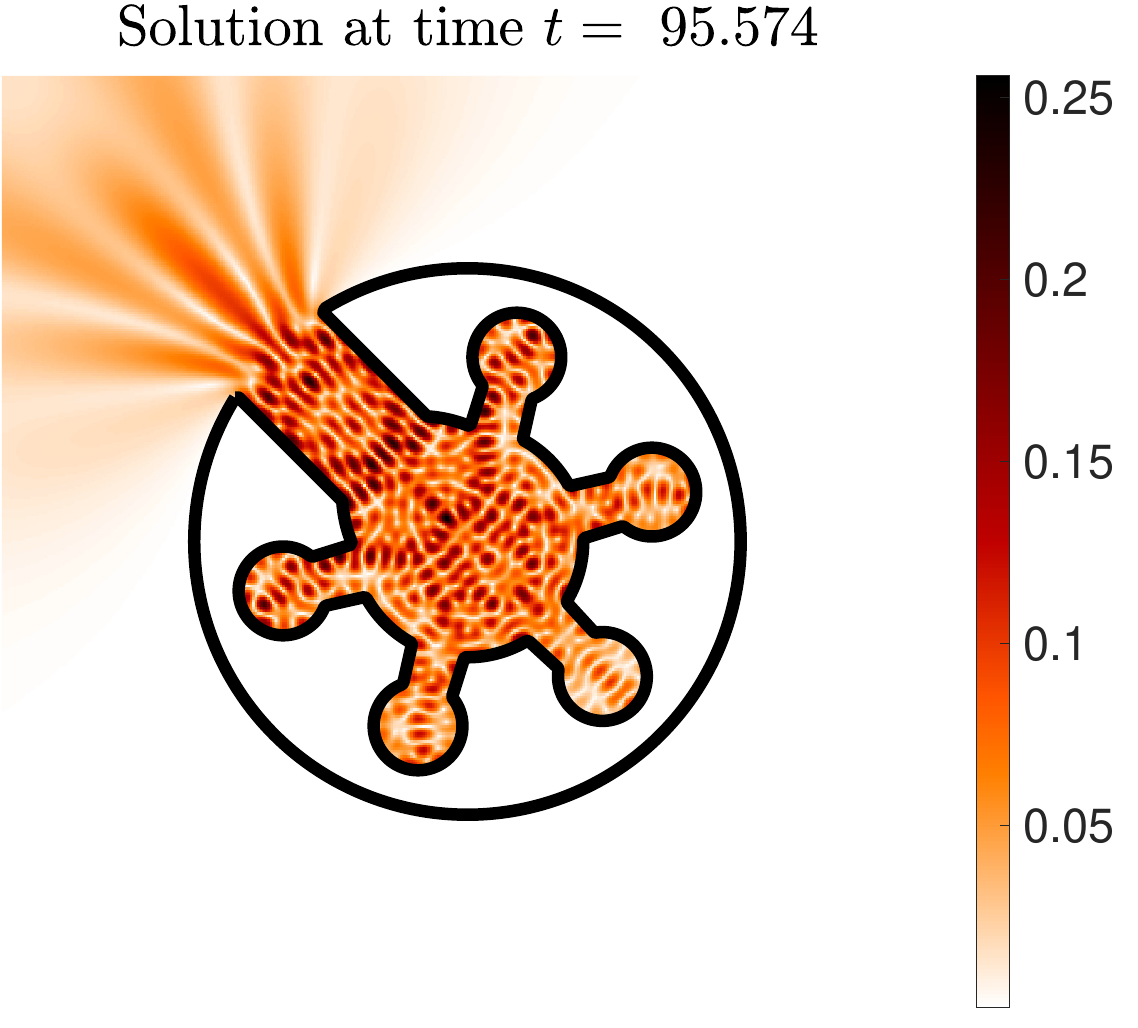}
\caption{ Magnitude of total field in a magnetron domain at $t = 31.311$ (left), $t = 42.787$ (center) and $t = 95.574$ (right)}
\label{fig:cornersrad}
 \end{figure}

We consider two domains: the absolute value of the total field is shown in a keyhole domain in \Cref{fig:cornerskey}, and in a more challenging magnetron domain in \Cref{fig:cornersrad}.  A different magnetron domain with thinner channels is shown in \Cref{fig:gasket_magplots}. For these experiments, $\delta = .025$ and $m = 1000$.

\subsubsection{Implementation details for nonsmooth domains}

To handle corners, we apply a standard refinement scheme to the discretization of the integral equation~\eqref{eq:mixedBIE} and the contour integral~\eqref{eq:U_int} that defines $\hat{U}(x, \omega)$. As discussed in~\cite{martinsson2019fast}, the scatterer is divided into several pieces, called panels, and Gauss-Legendre quadrature is applied to each panel with Kolm-Rokhlin correction~\cite{kolm2001numerical} applied to manage singularities in the layer potentials. To improve accuracy, the panels are refined geometrically into the corners. This basic method is sufficient for general 2D domains and can also be adapted to smooth domains with regions of high curvature. Alternative methods include QBX-type quadratures~\cite{klockner2013quadrature} and other specialized corner quadrature schemes~\cite{bruno2009high,hoskins2019numerical, serkh2019solution}.    

\subsection{A 3D example}
\begin{figure}[h!]
\centering
\includegraphics[scale = .4,trim={1cm 4cm 1cm 3cm},clip=true]{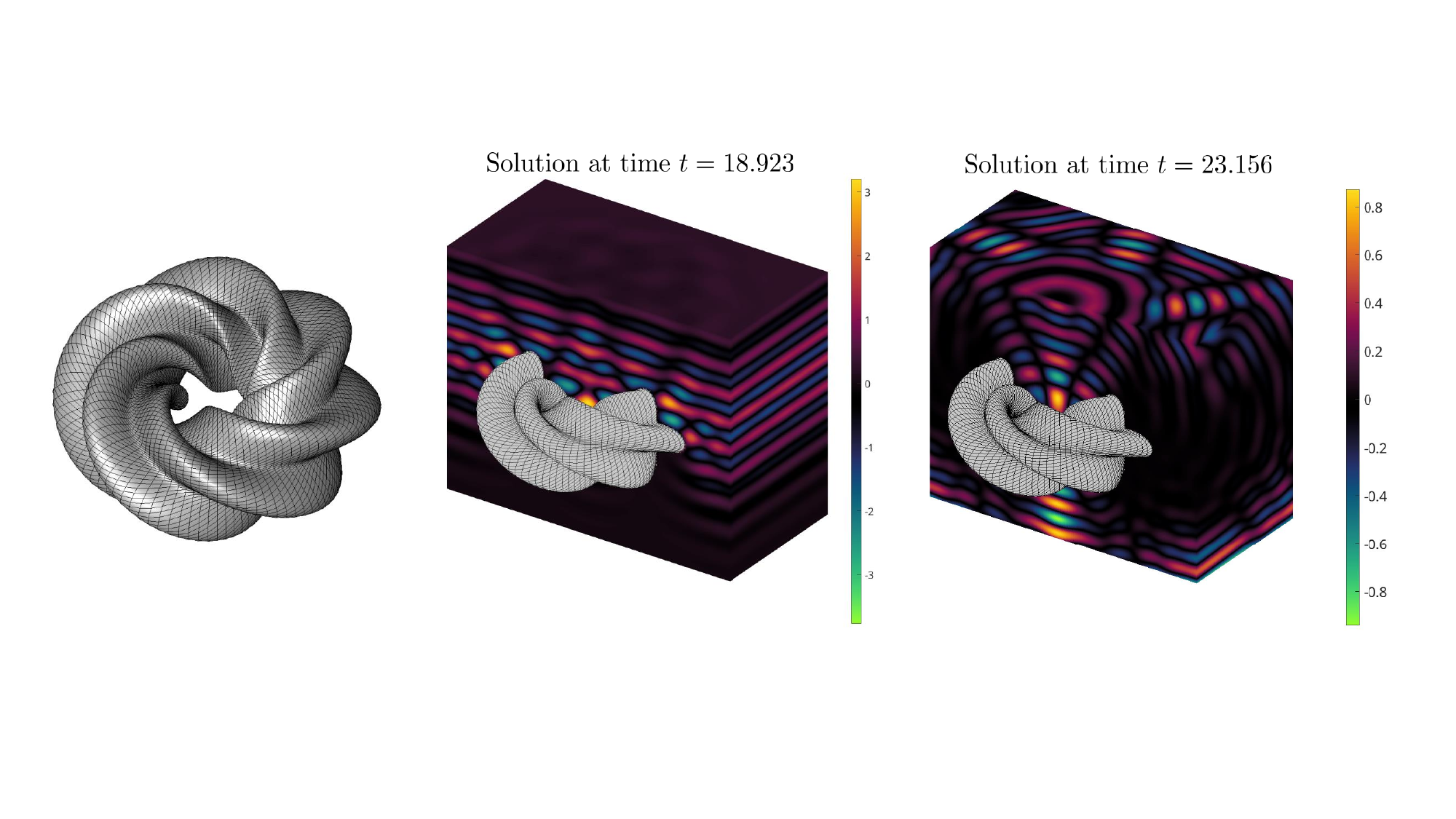}
\caption{A biperiodic 3D ``cruller" scatterer (left), and the real part of the total field at times $t = 18.923$ (center) and $t = 23.156 $ (right).}
\label{fig:cruellerplots}
\end{figure}

Since the damping+correction scheme only involves contour integration with respect to the frequency parameter $\omega$, it can be applied straightforwardly when $\Omega$ is a 3D domain. Nothing about the inverse Fourier transform changes, but it is substantially more computationally intensive to evaluate $\hat{U}(x, \omega)$.   As in the 2D setting, boundary integral methods can be used to solve for each density function associated with a given frequency, so that with the FMM, one can quickly evaluate $\hat{U}(x_k, \omega_j)$. We use the FMMLIB3D package~\cite{greengard2020fmmlib3d}. Specialized fast quadrature schemes for surface integrals  are needed to work with such domains. Our example in \Cref{fig:cruellerplots} is a proof-of-concept that uses a biperiodic ``cruller" domain from the BIE3D package~\cite{BIE3d} and Zeta-correction quadrature~\cite{wu2021corrected}. In future work, we plan to implement 3D examples involving more sophisticated domains with severe trapping.

\section{Acknowledgments} We thank Vladimir Rokhlin for multiple helpful conversations and for pointing out his early work using contour deformation to solve scattering problems in~\cite{rokhlin1983solution}. We thank Adrianna Gillman for sharing an implementation of Kolm-Rokhlin quadrature~\cite{kolm2001numerical} for nonsmooth contours, and we thank Alex Barnett for many helpful suggestions and conversations. 

\bibliographystyle{siamplain}
\bibliography{refs}

\appendix

\section{Code for scatterers} \label{apx:scatter}
MATLAB code for producing each of the scatterers used in our examples is supplied below. We use Chebfun~\cite{Chebfun} because it nicely handles the computation derivatives and normal directions, as well as plotting, in subsequent stages. All code used in our experiments, including a gallery command for scatterers, can be found at~\cite{TFAW}. 

\begin{enumerate}
   
    \item smooth C-curves (as in \Cref{fig:c_curveplots}): 
     \begin{verbatim}
a = 3;      % bigger = sharper corners
b = 2.8;    % smaller = narrower cavity opening
c = .1;     % larger = fatter all around
d = 1;      % rotation parameter

t = chebfun('t', [0, 2*pi]);
r = 3+c*tanh(a*cos(t));
th = b*sin(t);

gam = d*exp(1i*th).*r;

    \end{verbatim}
\item Crescents (as in \Cref{fig:crescent}):
    \begin{verbatim}
r = 5;              % size param
th = exp(1i*pi/2);  % angle param
p = 1i;            % relative shift
a1 = .24; a2 = .9;  % crescent shape params
d = 3;             % controls distance between crescents
    
t = chebfun('t', [0, pi],'splitting','on');
cres = @(t) exp(-2i*t)- a1/(exp(-2i*t)+a2) + d/2; % displaced crescent

gam = cres(t);
gam = join(gam+p, -cres(pi+t));
gam = r*th*gam;

    \end{verbatim} 
\item Keyhole domain (as in \Cref{fig:poleplot} and \Cref{fig:cornerskey}): 
    \begin{verbatim}
r = 2;      % inner rad
R = 3;      % outer rad
e = .3;     % thickness param
theta = pi; % controls rotation
c = [-R+e*1i -r+e*1i -r-e*1i -R-e*1i];

t = chebfun('t',[0 2*pi/4],'splitting','on');
l = @(t,a,b) a + (t/(2*pi)*4)*(b-a); % line segment
circ = @(t,a,b) a.^(1 - t/(2*pi)*4) .* b.^(t/(2*pi)*4); % circle segment

gam = join(l(t,c(1),c(2)),circ(t,c(2),c(3)),l(t,c(3),c(4)),circ(t,c(4),c(1)));
gam = exp(1i*(pi+theta)) * gam;

    \end{verbatim}      
\item magnetron domain (as in \Cref{fig:gasket_magplots} and bottom row of \Cref{fig:cornersrad}):
    \begin{verbatim}
% Define parameters
r = 1.5;            % radii small circs.
R = 5;              % radius large outer circ.
e = .8;
outr = .8;
deltatheta = .15;
n = 6-1;            % n-1 = number of sm. circs. 
Btheta = 5*pi/4;    % controls rotations

% Define inner smaller circle centers
thetastilde = 2 * (0:n) * pi / (n+1); 
thetas = [thetastilde(thetastilde >= pi),thetastilde(thetastilde < pi)];

centers = (R + r) / 2 * exp(1i * thetas); 

% Remove center point at -(R+r)/2
thetas(real(centers) == -(R + r)/2) = [];
centers(real(centers) == -(R + r)/2) = [];
n = length(centers);

% Define main keyhole points
c = [-R - e*1i, -r - e*1i];
for i1 = 1 : n
    theta = thetas(i1);
    cent = centers(i1);
    c = [c,r*exp((theta-deltatheta)*1i),cent +...
    outr * exp((theta+r/outr*deltatheta+pi)*1i),cent +...
    outr * exp((theta-r/outr*deltatheta+pi)*1i),...
    r*exp((theta+deltatheta)*1i)];
end

c = [c, -r + e*1i, -R + e*1i];
NN = length(c);
a = 2*pi/NN;

% Define parameter for curve
t = chebfun('t', [0, a], 'splitting', 'on');
l = @(t,c1,c2) c1 + (t/a)*(c2-c1); % line segment
circ1 = @(t,c1,c2) c1.^(1 - t/a) .* c2.^(t/a); % circle segment
circ2 = @(t,c,th) c + outr * exp(1i * ((1 - t/a) * th(1) + (t/a) * th(2)));

gam1 = l(t,c(1),c(2));
for i1 = 2:4:(length(c)-4)
    gam1 = join(gam1, circ1(t,c(i1),c(i1+1)));
    gam1 = join(gam1, l(t,c(i1+1),c(i1+2)));

    cent = centers(ceil(i1/4));
    angle_start = angle(c(i1+2) - cent);
    angle_end = angle(c(i1+3) - cent);
    if abs(angle_end-angle_start)<pi
        angle_end = angle_end + 2*pi;
    end
    angles = [angle_start,angle_end];
    
    gam1 = join(gam1, circ2(t,cent,angles));
    gam1 = join(gam1, l(t,c(i1+3),c(i1+4)));

end

gam1 = join(gam1, circ1(t,c(end-2),c(end-1)));
gam1 = join(gam1, l(t,c(end-1),c(end)));
gam1 = join(gam1, circ1(t,c(end),c(1)));

tfull = chebfun('t', [0, 2*pi], 'splitting', 'on');
gam = exp(1i*(pi+Btheta))*gam1(2*pi-tfull);
\end{verbatim}
\end{enumerate}

\end{document}